\documentclass[11pt,a4paper]{article}

\usepackage[T1]{fontenc}
\usepackage[utf8]{inputenc}
\usepackage{lmodern}
\usepackage{microtype}
\usepackage[margin=2.7cm]{geometry}
\usepackage{amsmath,amssymb,amsthm,url}
\usepackage{amscd}
\usepackage[all]{xy}
\usepackage{tikz-cd}
\usepackage{mathtools}
\usepackage[hidelinks]{hyperref}
\usepackage{bookmark}

\hypersetup{
  pdftitle={A Stone duality theorem for arbitrary lattices},
  pdfauthor={Andr\'es R\'ios},
  pdfsubject={Preprint},
  pdfkeywords={Lattice theory, duality, Stone representation theorem, bitopological spaces, Balbes--Dwinger spaces}
}

\numberwithin{equation}{section}

\theoremstyle{plain}
\newtheorem{theorem}{Theorem}[section]
\newtheorem{lemma}[theorem]{Lemma}
\newtheorem{proposition}[theorem]{Proposition}
\newtheorem{corollary}[theorem]{Corollary}

\theoremstyle{definition}
\newtheorem{definition}[theorem]{Definition}
\newtheorem{notation}[theorem]{Notation}
\newtheorem{remark}[theorem]{Remark}
\newtheorem{example}[theorem]{Example}

\newenvironment{Th}{\begin{theorem}}{\end{theorem}}
\newenvironment{Lem}{\begin{lemma}}{\end{lemma}}
\newenvironment{Prop}{\begin{proposition}}{\end{proposition}}
\newenvironment{Cor}{\begin{corollary}}{\end{corollary}}
\newenvironment{Def}{\begin{definition}}{\end{definition}}

\title{\textbf{A Stone duality theorem\\for arbitrary lattices}}
\author{Andr\'es R\'ios\\[0.4em]
  {\small Departamento de Matem\'aticas, Universidad El Bosque}\\
  {\small Bogot\'a, Colombia}\\
  {\small\href{mailto:afriosm@unbosque.edu.co}{\texttt{afriosm@unbosque.edu.co}}}}
\date{\small Preprint --- September 12, 2026}

\begin{document}
\maketitle

\begin{abstract}
We establish a Stone duality theorem for arbitrary lattices, without assuming
distributivity or boundedness. Rather than extending Priestley-style ordered
representations, our construction follows the modern Stone--Balbes--Dwinger
viewpoint: a lattice is recovered from distinguished open subsets of a
spectrum. Since prime ideals are no longer sufficient in the non-distributive
setting, we replace them by comaximal pairs and associate to each lattice
\(L\) a bitopological spectrum \(\operatorname{spec}_B(L)\).

The two topologies arise from two natural Stone-like maps: one encodes the
join-semilattice structure of \(L\), while the other encodes its
meet-semilattice structure. These two maps coincide precisely in the
distributive case, showing that the single topology of the classical
Stone--Balbes--Dwinger representation is the collapsed form of a genuinely
bitopological construction. We introduce pairwise Balbes--Dwinger spaces as
the corresponding class of bitopological spaces and prove that \(L\) is
naturally isomorphic to the lattice of essential subsets of
\(\operatorname{spec}_B(L)\). This yields a representation theorem for
arbitrary lattices. We then study functoriality and prove a dual equivalence
between lattices with quasi-proper homomorphisms and pairwise
Balbes--Dwinger spaces with the appropriate morphisms. In the distributive
case, the two topologies coincide and the classical Stone--Balbes--Dwinger
duality is recovered.
\end{abstract}



\clearpage
\tableofcontents
\clearpage

\section{Introduction}

In 1938, Marshall H.~Stone \cite{Stone} initiated a systematic program for
representing distributive lattices by means of topological spaces. In his
seminal work, he established topological representation theorems for
distributive lattices and Brouwerian logics, laying the foundations for what
is now known as Stone-type duality. Later, Balbes and Dwinger
\cite{Balbes-Dwinger} extended Stone's ideas to the arbitrary, not necessarily
bounded, distributive case. They introduced a class of topological spaces,
which we call \emph{Balbes--Dwinger spaces}, and showed how a distributive
lattice can be recovered from distinguished open subsets of its spectrum.
Moreover, they described how algebraic properties of lattices translate into
topological properties and conversely.

However, the Balbes--Dwinger representation is not fully functorial: not
every homomorphism of distributive lattices induces a continuous map between
the corresponding spaces. To address this issue, Acosta \cite{Acosta}
restricted the morphisms in the category of distributive lattices by
introducing the notion of a \emph{proper lattice homomorphism}, and
derived---building on Balbes and Dwinger---a coequivalence between the category
of distributive lattices with proper homomorphisms and the category of
Balbes--Dwinger spaces.

A different path was opened by Priestley
\cite{Priestley70,Priestley72}, who represented bounded distributive lattices
by ordered topological spaces, now called \emph{Priestley spaces}. This
ordered representation has the advantage of yielding a full duality for
bounded distributive lattices. Urquhart \cite{Urquhart} later extended
Priestley's representation theory to the non-distributive setting by
introducing doubly ordered topological spaces, which we call
\emph{Urquhart spaces}. His construction provides a representation theory for
bounded lattices, but it also fails to be fully functorial unless morphisms
are suitably restricted. Mancera and Acosta \cite{Mancera} resolved this
difficulty by restricting morphisms in the category of bounded lattices,
obtaining a duality theorem based on Urquhart's ideas. Another approach to
recovering functoriality was developed by Allwein and Hartonas
\cite{AllHart}, who enlarged the underlying spaces by considering additional
points, obtaining a full duality for bounded lattices.

The aim of the present paper is different from the Priestley--Urquhart line.
We do not begin from ordered representations. Instead, we return to the
Stone--Balbes--Dwinger viewpoint: a lattice should be recovered from suitable
open subsets of a spectrum. The main obstruction is that prime ideals are no
longer adequate in the non-distributive case; indeed, non-distributive
lattices may have no prime ideals at all. We therefore replace prime ideals
by \emph{comaximal pairs} of ideals and filters, following an idea already
present in Urquhart's work, and we use them as the points of a new spectrum.

The crucial feature of the construction is that, for an arbitrary lattice,
there are two natural Stone-like maps. One of them encodes the
join-semilattice structure, while the other encodes the meet-semilattice
structure. In the distributive case these two maps coincide, and this
coincidence characterizes distributivity. Thus the single topology appearing
in the classical Stone--Balbes--Dwinger representation is not lost in the
general case; rather, it splits naturally into two interacting topologies.
This leads us to associate to each lattice \(L\) a bitopological spectrum
\[
\operatorname{spec}_B(L),
\]
whose underlying set is the set of comaximal pairs of \(L\).

To characterize the spaces that arise in this way, we introduce
\emph{pairwise Balbes--Dwinger spaces}. These are bitopological spaces
designed so that the two topologies recover, respectively, the join and meet
behaviour of a lattice, while the transition maps between the two topologies
glue these pieces back into a genuine lattice structure. For every lattice
\(L\), the spectrum \(\operatorname{spec}_B(L)\) is a pairwise
Balbes--Dwinger space, and the lattice of essential subsets of
\(\operatorname{spec}_B(L)\) is naturally isomorphic to \(L\). This yields a
representation theorem for arbitrary lattices, without assuming
distributivity or boundedness.

As in the Balbes--Dwinger and Urquhart settings, full functoriality cannot be
obtained with arbitrary lattice homomorphisms. Following the strategy of
Acosta \cite{Acosta} and Mancera \cite{Mancera}, we introduce a suitable
restriction on morphisms, called \emph{quasi-proper homomorphisms}. With this
choice, the spectrum construction becomes functorial and yields a
coequivalence between the category of arbitrary lattices with quasi-proper
homomorphisms and the category of pairwise Balbes--Dwinger spaces with the
appropriate strongly bicontinuous morphisms.

Finally, we show that the distributive case is exactly the case in which the
bitopological construction collapses to an ordinary topological one. More
precisely, distributive lattices correspond to those pairwise
Balbes--Dwinger spaces whose two topologies coincide. These are precisely the
doubly Balbes--Dwinger spaces, and under this identification the classical
Stone--Balbes--Dwinger duality is recovered. Thus the theory developed here
is not merely analogous to Stone duality: it extends the
Stone--Balbes--Dwinger representation to arbitrary lattices, with the
classical distributive theory appearing as the special case in which the two
topologies collapse into one.

\section{Preliminary results}
\subsection{The Stone duality theorem}
In this section we briefly introduce the main ideas behind the Stone duality theorem. It is important to note that the version presented here is not the classical formulation given by Stone himself; rather, we follow the modern approach of \cite{Balbes-Dwinger} and \cite{Acosta}, which nonetheless implies the classical Stone duality of \cite{Stone}.

Let $L$ be a distributive lattice. Following \cite{Acosta} and \cite{Balbes-Dwinger}, the idea behind the Stone duality theorem is to associate to $L$ a topological space, namely the \textbf{spectrum of $L$}, whose points are the prime ideals of $L$; that is, the underlying set of this space is $spec(L)$. To define a topology on $spec(L)$ we consider the function
\[
d_L : L \to \wp(spec(L)), \qquad 
d_L(x) = \{\, P \in spec(L) \mid x \notin P \,\}.
\]
The image of $d_L$ is closed under finite intersections, hence $Im(d_L)$ is a basis for a topology on $spec(L)$, called the \emph{Zariski topology}. We continue to denote this topological space by $spec(L)$, without further mention of the topology.

It is important to observe that if $L$ is a non-trivial distributive lattice, then $spec(L)$ is non-empty. This follows from the prime ideal theorem for distributive lattices \cite{Acosta}. However, the same is not true in the non-distributive case: there exist non-distributive lattices, such as the five-element diamond \(\mathfrak{M}_5\), whose spectrum is empty.

The ``topologizing'' function $d_L$ is an injective lattice homomorphism, and therefore induces an isomorphism between $L$ and $Im(d_L)$. Thus, there is an isomorphic copy of $L$ inside the topology of $spec(L)$. The first step in Stone's representation theory is to recover this copy by identifying the topological properties that characterize $Im(d_L)$ inside the lattice of open subsets of $spec(L)$. This is achieved by showing that, for any $x \in L$, the set $d_L(x)$ is compact and open. The openness is immediate from the definition of the topology; compactness follows from the Balbes–Dwinger theorem (stated below). Hence a natural idea is to recover $L$ from the collection of compact-open subsets of $spec(L)$. A difficulty arises, however: the empty set is always compact-open, and thus the ordered set of compact-open subsets always has a minimum element, even if $L$ does not. This motivates the following definition.

\begin{Def}[\cite{Acosta}]
Let $X$ be a topological space. A non-empty subset of $X$ is \textbf{fundamental} if it is compact and open. The empty set is \textbf{fundamental} if every family of compact-open subsets of $X$ with the finite intersection property has non-empty intersection. We denote by $\mathfrak{F}(X)$ the collection of fundamental subsets of $X$.
\end{Def}

It can be shown that $L$ has a minimum element (i.e.\ $0$) if and only if the empty set is a fundamental subset of $spec(L)$. Moreover, a key fact about the spectrum of a distributive lattice is that every fundamental subset of $spec(L)$ is of the form $d_L(x)$ for a unique $x \in L$. Thus the image of $d_L$ coincides with $\mathfrak{F}(spec(L))$, and $d_L$ induces an isomorphism between $L$ and $\mathfrak{F}(spec(L))$. This reconstruction of $L$ from the fundamental subsets of its spectrum is usually described as the Stone representation theorem, as it shows that any distributive lattice is a sublattice of $\wp(X)$ for some set $X$ (namely $X = spec(L)$).

So far we have associated to each distributive lattice a topological space and described how the original lattice can be recovered from the topology. The next step in Stone duality is to understand this construction functorially. An obstruction arises: in general, an arbitrary lattice homomorphism does not preserve prime ideals by inverse image (though this holds in certain special settings, such as bounded distributive lattices with bounded homomorphisms). Following \cite{Acosta}, one resolves this by restricting the allowed morphisms. A lattice homomorphism $f : L \to N$ is called \emph{proper} if for every prime ideal $P \subseteq N$, the preimage $f^{-1}(P)$ is a prime ideal of $L$. Clearly, any proper homomorphism induces a function
\[
f^{-1} : spec(N) \to spec(L),
\]
which can be shown to be \emph{strongly continuous} (i.e.\ continuous and with the property that the preimage of any fundamental subset of $spec(L)$ is fundamental in $spec(N)$). Thus, if $\mathfrak{D}_p$ is the category of distributive lattices with proper homomorphisms, the spectrum construction defines a functor
\[
spec : \mathfrak{D}_p \to Top,
\]
where $Top$ is the category of topological spaces.

The next step is to characterize the essential image of this functor; that is, to identify the topological spaces that arise as spectra of distributive lattices. The key properties of $spec(L)$ are the following (see \cite{Balbes-Dwinger, Acosta}):
\begin{enumerate}
    \item[(i)] $spec(L)$ is a $T_0$ space.
    \item[(ii)] $spec(L)$ is coherent (the fundamental subsets form a basis closed under finite intersections).
    \item[(iii)] The fundamental subsets of $spec(L)$ are bireducible (arbitrary intersections contained in arbitrary unions reduce to finite ones).
\end{enumerate}

The third property is the least familiar. Its explicit formulation for $spec(L)$ is given by the Balbes–Dwinger theorem:

\begin{Th}[\textbf{Balbes--Dwinger} {\cite{Acosta}}]
Let $L$ be a distributive lattice. If $V, W \subseteq L$ are non-empty sets such that
\[
\bigcap_{x \in W} d_L(x) \subseteq \bigcup_{y \in V} d_L(y),
\]
then there exist finite subsets $V_0 \subseteq V$ and $W_0 \subseteq W$ such that
\[
\bigcap_{x \in W_0} d_L(x) \subseteq \bigcup_{y \in V_0} d_L(y).
\]
\end{Th}

The name of this theorem is not standard in the literature; we follow the terminology of \cite{Acosta}.

A topological space $X$ is a \emph{Balbes--Dwinger space} if it is $T_0$, coherent, and its fundamental subsets are bireducible. Although this terminology is also non-standard, it is clear that the spectrum of any distributive lattice is a Balbes--Dwinger space. If $X$ is such a space, coherence implies that $\mathfrak{F}(X)$ is a distributive lattice (with intersection and union as the lattice operations). Hence we may consider its spectrum $spec(\mathfrak{F}(X))$, which is again a Balbes--Dwinger space. To relate $X$ with $spec(\mathfrak{F}(X))$, notice that for any $x \in X$ the set
\[
I_x = \{\, A \in \mathfrak{F}(X) \mid x \notin A \,\}
\]
is a prime ideal of $\mathfrak{F}(X)$ \cite{Acosta}. We may therefore define a map
\[
h_X : X \to spec(\mathfrak{F}(X)), \qquad h_X(x) = I_x,
\]
which is a homeomorphism \cite{Acosta}. Thus the Balbes--Dwinger spaces are precisely the objects in the essential image of the functor $spec$. As before, to understand the assignment $X \mapsto \mathfrak{F}(X)$ functorially, we must restrict morphisms: continuous maps alone are not sufficient. Let $\mathfrak{BD}$ denote the category whose objects are Balbes--Dwinger spaces and whose morphisms are strongly continuous maps. Any such map $f : X \to Y$ induces (by inverse image) a function
\[
f^{*} : \mathfrak{F}(Y) \to \mathfrak{F}(X),
\]
which is a proper lattice homomorphism. Thus,
\[
\mathfrak{F} : \mathfrak{BD} \to \mathfrak{D}_p
\]
is a functor.

We can now state the Stone duality theorem in its categorical form.

\begin{Th}[\textbf{Stone duality theorem}]\label{StoneDual}
The functor $spec : \mathfrak{D}_p \to \mathfrak{BD}$ is a coequivalence of categories, whose dual functor is $\mathfrak{F} : \mathfrak{BD} \to \mathfrak{D}_p$.
\end{Th}

The main aim of this paper is to generalize this theorem to the setting of arbitrary lattices, including the non-distributive case. This cannot be achieved using the techniques presented above, since the spectrum of a non-distributive lattice may be empty. We must therefore develop new methods to overcome this obstruction and obtain a duality in full generality.

\subsection{Priestley and Urquhart representation theorems}
In this section we present the main ideas behind the Priestley and Urquhart representation theorems for lattices. We will not emphasize technical details of these constructions; rather, we briefly outline their conceptual foundations, since they serve only as inspiration for our work. The main focus of this text is to extend the techniques and ideas of the Stone-type duality for distributive lattices developed in the previous section.

First, in Priestley's representation theory for distributive lattices \cite{Priestley70,Priestley72}, the idea is to associate to any bounded distributive lattice $L$ an ordered topological space whose underlying set is the set of prime ideals of $L$, denoted $spec_P(L)$. As in the Stone duality framework, one can recover an isomorphic copy of $L$ inside the topology of $spec_P(L)$ by combining order-theoretic and topological properties of this space. This construction yields a functor from the category of bounded distributive lattices (with bounded lattice homomorphisms) to the category of ordered topological spaces, whose essential image is the category of Priestley spaces---that is, compact ordered topological spaces satisfying the Priestley separation axiom \cite{Priestley70,Priestley72}. As in Stone duality, this leads to a dual equivalence between the category of bounded distributive lattices and the category of Priestley spaces (with continuous order-preserving maps). A key advantage of Priestley duality is that it represents \emph{all} morphisms between bounded distributive lattices; no restriction on lattice morphisms is required.

Urquhart's work generalizes Priestley's theory to the case of bounded, not necessarily distributive, lattices. To each bounded lattice $L$, one associates a bitopological space equipped with two partial orders, called the \emph{Urquhart spectrum} of $L$ and denoted $spec_U(L)$. As in previous representation theories, one can characterize those doubly ordered topological spaces that arise as the Urquhart spectrum of some bounded lattice. One may also study this representation functorially by choosing appropriate categories of bounded lattices and appropriate categories of doubly ordered topological spaces (after imposing certain restrictions on morphisms). For detailed treatments we refer to \cite{Mancera} and \cite{Urquhart}.

Since there is a natural way to pass from the Priestley representation to the Stone representation in the setting of bounded lattices, Urquhart's work suggests that one should expect an extension of Stone duality using purely topological methods (that is, without relying on order structure on the topological side). This observation is the starting point of the present work.

\section{The construction of the spectrum}
In this section we present the construction of the spectrum of an arbitrary lattice $L$. We begin by defining the underlying set of the space. Then we introduce suitable topologies on this set and study some of their basic properties. For this purpose, we first define two essential functions that allow us to ``bitopologize'' the underlying set. We conclude the section by showing how any lattice can be reconstructed from its spectrum, thereby establishing a representation theorem for arbitrary lattices.

\subsection{The underlying set}
Following \cite{Cornish}, there exists an equivalence between the Stone and Priestley representations for bounded distributive lattices. This equivalence preserves underlying sets; that is, for a bounded distributive lattice, the underlying set of its Stone spectrum and the underlying set of its Priestley spectrum coincide: both are the set of prime ideals of the lattice. Since Urquhart established an extension of Priestley duality to the setting of bounded, not necessarily distributive, lattices \cite{Urquhart}, the equivalence between the Priestley and Stone representations suggests that the underlying set of any hypothetical extension of Stone's construction to arbitrary lattices should coincide with the underlying set of its Urquhart spectrum.

Let $L$ be a lattice.

\begin{Def}[\cite{Urquhart}, Definition~2.1]
A pair $(I,F)$ consisting of an ideal $I$ and a filter $F$ is a \textbf{comaximal pair} of $L$ if:
\begin{enumerate}
    \item[(i)] $I \cap F = \emptyset$,
    \item[(ii)] if $J$ is an ideal of $L$ with $I \subsetneq J$, then $J \cap F \neq \emptyset$,
    \item[(iii)] if $K$ is a filter of $L$ with $F \subsetneq K$, then $I \cap K \neq \emptyset$.
\end{enumerate}
The set of comaximal pairs of $L$ is denoted by $\mathfrak{M}(L)$.
\end{Def}

The underlying set of the Urquhart spectrum is the set $\mathfrak{M}(L)$ of comaximal pairs. As mentioned above, this will also be the underlying set of the spectrum introduced in this text.

Any prime ideal $P$ of $L$ determines a comaximal pair $(P, P^c)$. This defines a function
\[
b_L : spec(L) \longrightarrow \mathfrak{M}(L), \qquad 
P \longmapsto (P, P^c),
\]
which is always injective (provided that $spec(L)$ is non-empty). In fact, $b_L$ provides a criterion for determining when a lattice is distributive.

The following result captures the fact that there are enough comaximal pairs
to separate disjoint ideals and filters. It may be viewed as a non-distributive
analogue of the prime ideal theorem.

\begin{Prop}[\textbf{Comaximal pair theorem} \cite{Urquhart}]\label{Com}
Let $L$ be any lattice. If $I$ is an ideal of $L$ and $F$ is a filter of $L$ with $I \cap F = \emptyset$, then there exist an ideal $J \supseteq I$ and a filter $K \supseteq F$ such that $(J, K)$ is a comaximal pair of $L$.
\end{Prop}

A key consequence is that if $L$ has more than one element, then $\mathfrak{M}(L) \neq \emptyset$. Whenever $L$ has at least two distinct elements, one can always find an ideal and a filter that are disjoint; then Proposition~\ref{Com} yields a comaximal pair. This stands in sharp contrast with the behavior of prime ideals: there exist non-distributive lattices with no prime ideals at all, such as the
five-element diamond \(\mathfrak{M}_5\). Thus comaximal pairs extend the notion of prime ideals to the general, possibly non-distributive, setting, preserve their separation properties, and always exist in any lattice with at least two points, making them particularly well suited for a general representation theory.

\begin{Lem}\label{ComaximalPrimeCriterion}
A lattice \(L\) is distributive if and only if every ideal which occurs as
the first component of a comaximal pair of \(L\) is prime.
\end{Lem}

\begin{proof}
Suppose first that \(L\) is distributive, and let \((I,F)\) be a comaximal
pair of \(L\). By the prime ideal theorem for distributive lattices, there
exists a prime ideal \(P\) of \(L\) such that
\[
    I\subseteq P
    \qquad\text{and}\qquad
    P\cap F=\emptyset.
\]
By the maximality condition on \(I\), we must have \(I=P\). Hence \(I\) is
prime.

Conversely, suppose that every ideal which occurs as the first component of
a comaximal pair of \(L\) is prime. We prove that \(L\) is distributive.
Assume, toward a contradiction, that \(L\) is not distributive. By Birkhoff's characterization of distributive lattices
\cite[Chapter~II, Section~1]{Gratzer}, \(L\) contains a sublattice
isomorphic to either the five-element diamond \(\mathfrak{M}_5\) (often denoted by \(M_3\))
or the five-element pentagon \(\mathfrak{N}_5\).

First suppose that \(L\) contains a sublattice isomorphic to
\(\mathfrak{M}_5\). Write its elements as
\[
    0_S,\ a,\ b,\ c,\ 1_S,
\]
where \(a,b,c\) are the three pairwise incomparable middle elements. Let
\[
    I=(a]
    \qquad\text{and}\qquad
    F=[b)
\]
be the principal ideal and the principal filter generated in \(L\). Since \(a\) and \(b\) are incomparable, the principal ideal \((a]\) and the
principal filter \([b)\) are disjoint: otherwise there would be an element
\(u\in L\) such that
\[
    b\leq u\leq a,
\]
which would imply \(b\leq a\), a contradiction. By
Proposition~\ref{Com}, there exists a comaximal pair \((Q,G)\) of \(L\) such
that
\[
    I\subseteq Q
    \qquad\text{and}\qquad
    F\subseteq G.
\]
We claim that \(Q\) is not prime. Indeed, since \(0_S=b\wedge c\) and
\(0_S\leq a\), we have
\[
    b\wedge c=0_S\in Q.
\]
Moreover, \(b\notin Q\), because \(b\in F\subseteq G\) and \(Q\cap G=\emptyset\).
Also \(c\notin Q\). For if \(c\in Q\), then
\[
    1_S=a\vee c\in Q,
\]
and since \(b\leq 1_S\), the ideal property would imply \(b\in Q\), a
contradiction. Thus \(b\wedge c\in Q\), while \(b\notin Q\) and \(c\notin Q\).
Therefore \(Q\) is not prime, contradicting the hypothesis.

Now suppose that \(L\) contains a sublattice isomorphic to
\(\mathfrak{N}_5\). Write its elements as
\[
    0_S,\ a,\ b,\ c,\ 1_S,
\]
with \(c<b<1_S\), \(0_S<a<1_S\), and
\[
    a\wedge b=0_S,
    \qquad
    a\vee c=1_S.
\]
Let
\[
    I=(c]
    \qquad\text{and}\qquad
    F=[b).
\]
Since \(c<b\), the principal ideal \((c]\) and the principal filter \([b)\)
are disjoint: otherwise there would be an element \(u\in L\) such that
\(b\leq u\leq c\), which would imply \(b\leq c\), a contradiction. By Proposition~\ref{Com},
there exists a comaximal pair \((Q,G)\) of \(L\) such that
\[
    I\subseteq Q
    \qquad\text{and}\qquad
    F\subseteq G.
\]
Again we show that \(Q\) is not prime. Since \(0_S=a\wedge b\) and
\(0_S\leq c\), we have
\[
    a\wedge b=0_S\in Q.
\]
Moreover, \(b\notin Q\), because \(b\in F\subseteq G\) and
\(Q\cap G=\emptyset\). Also \(a\notin Q\). Indeed, if \(a\in Q\), then
\[
    1_S=a\vee c\in Q,
\]
and since \(b\leq 1_S\), the ideal property would imply \(b\in Q\), a
contradiction. Hence \(a\wedge b\in Q\), while \(a\notin Q\) and \(b\notin Q\).
Thus \(Q\) is not prime, again contradicting the hypothesis.

In both cases we obtain a contradiction. Therefore \(L\) must be distributive.
\end{proof}

\begin{Prop}\label{PCLat}
The lattice \(L\) is distributive if and only if \(b_L\) is bijective; that is,
if and only if every comaximal pair of \(L\) is of the form \((P,P^c)\) for
some prime ideal \(P\) of \(L\).
\end{Prop}

\begin{proof}
The map \(b_L\) is injective by definition. Indeed, if
\[
    b_L(P)=b_L(Q),
\]
then \((P,P^c)=(Q,Q^c)\), and hence \(P=Q\).

Suppose first that \(L\) is distributive, and let \((I,F)\) be a comaximal
pair of \(L\). By Lemma~\ref{ComaximalPrimeCriterion}, the ideal \(I\) is
prime. Hence \(I^c\) is a filter disjoint from \(I\). Since \(I\cap F=\emptyset\),
we have \(F\subseteq I^c\). By the maximality condition on \(F\), it follows
that
\[
    F=I^c.
\]
Therefore
\[
    (I,F)=(I,I^c)=b_L(I),
\]
and so \(b_L\) is surjective. Hence \(b_L\) is bijective.

Conversely, suppose that \(b_L\) is bijective. We show that every ideal which
occurs as the first component of a comaximal pair is prime. Let \(I\) be such
an ideal, and choose a filter \(F\) such that \((I,F)\) is a comaximal pair
of \(L\). Since \(b_L\) is surjective, there exists a prime ideal \(P\) of
\(L\) such that
\[
    (I,F)=b_L(P)=(P,P^c).
\]
Thus \(I=P\), and therefore \(I\) is prime. By
Lemma~\ref{ComaximalPrimeCriterion}, \(L\) is distributive.
\end{proof}

Therefore, for any distributive lattice $L$, every comaximal pair arises from a prime ideal. Thus, for distributive $L$, one may identify $spec(L)$ with $\mathfrak{M}(L)$ as sets, showing that comaximal pairs generalize prime ideals in the distributive case. Moreover, comaximal pairs also generalize the behavior of separating ideals and filters in arbitrary lattices.

\subsection{Two essential functions}
According to \cite{Cornish}, the order of the Priestley spectrum of a distributive lattice is the topological order induced by the topology of its Stone spectrum (the topological order is defined by $x \leq y$ if $y \in \overline{\{x\}}$ in this topology). The Urquhart spectrum consists of a topological space equipped with two orders, and these ideas suggest that it is possible to define two topologies on $\mathfrak{M}(L)$ inducing these two orders. We follow this suggestion; thus the spectrum introduced here will be a bitopological space (in the sense of Kelly \cite{Kelly}). Mimicking the construction of the Stone spectrum presented in the first section, we define two functions that will allow us to introduce two topologies on the set of comaximal pairs of a lattice.

Let $L$ be a lattice. Define
\begin{eqnarray*}
    \delta_L : L &\longrightarrow& \wp(\mathfrak{M}(L)) \\
    x &\longmapsto& \{(I,F)\in\mathfrak{M}(L)\mid x\notin I\},
\end{eqnarray*}
and
\begin{eqnarray*}
    \epsilon_L : L &\longrightarrow& \wp(\mathfrak{M}(L)) \\
    x &\longmapsto& \{(I,F)\in\mathfrak{M}(L)\mid x\in F\}.
\end{eqnarray*}

Before showing how these two functions induce topologies on $\mathfrak{M}(L)$, we study some of their basic properties.

\begin{Prop}\label{Copies}
\begin{enumerate}
    \item[(i)] $\delta_L$ and $\epsilon_L$ are injective functions,
    \item[(ii)] $\delta_L$ is a sup-semilattice homomorphism, i.e.\ $\delta_L(x\vee y)=\delta_L(x)\cup\delta_L(y)$,
    \item[(iii)] $\epsilon_L$ is an inf-semilattice homomorphism, i.e.\ $\epsilon_L(x\wedge y)=\epsilon_L(x)\cap\epsilon_L(y)$.
\end{enumerate}
\end{Prop}

\begin{proof}
\begin{enumerate}
    \item[(i)] We show that $\delta_L$ is injective; the proof for $\epsilon_L$ is analogous.  
    If $x\neq y$, then $x\not\le y$ or $y\not\le x$. Without loss of generality assume $x\not\le y$. Then the ideal $(y]$ and the filter $[x)$ are disjoint. By Proposition~\ref{Com}, there exists a comaximal pair $(I,F)$ such that $(y]\subseteq I$ and $[x)\subseteq F$. Hence $(I,F)\in\delta_L(x)$ but $(I,F)\notin\delta_L(y)$.

    \item[(ii)] This follows directly from the definition of ideal.

    \item[(iii)] This follows from the definition of filter.
\end{enumerate}
\end{proof}

Before constructing topologies on $\mathfrak{M}(L)$ induced by these two functions, it is useful to examine their behavior in the distributive case. If $L$ is distributive, all comaximal pairs are of the form $(P,P^c)$ for some prime ideal $P$, and in this case it is clear that $\delta_L=\epsilon_L$. Furthermore, the bijection $b_L$ between $spec(L)$ and $\mathfrak{M}(L)$ fits into the commutative diagram
\[
\xymatrix{
   & & spec(L) \ar[dd]^{b_L} \\
  L \ar@/^/[urr]^{d_L} \ar@/_/[drr]_{\delta_L=\epsilon_L} & & \\
 & & \mathfrak{M}(L)
}
\]
i.e.\ $\delta_L$ and $\epsilon_L$ coincide for distributive $L$, and in essence these functions generalize the function $d_L$ introduced in the first section.

\begin{Prop}\label{TopDis}
\(L\) is distributive if and only if \(\delta_L=\epsilon_L\).
\end{Prop}

\begin{proof}
Suppose first that \(L\) is distributive. By Proposition~\ref{PCLat}, every
comaximal pair of \(L\) is of the form \((P,P^c)\), where \(P\) is a prime
ideal. Hence, for every \(x\in L\),
\[
(P,P^c)\in\delta_L(x)
\Longleftrightarrow x\notin P
\Longleftrightarrow x\in P^c
\Longleftrightarrow (P,P^c)\in\epsilon_L(x).
\]
Thus \(\delta_L=\epsilon_L\).

Conversely, suppose that \(\delta_L=\epsilon_L\). Let \((I,F)\) be a
comaximal pair of \(L\). We show that \(I\) is prime. Suppose
\(x\wedge y\in I\) and \(x,y\notin I\). Then
\[
(I,F)\in\delta_L(x)\cap\delta_L(y).
\]
Since \(\delta_L=\epsilon_L\), and since \(\epsilon_L\) preserves finite
meets, we obtain
\[
(I,F)\in\epsilon_L(x)\cap\epsilon_L(y)
       =\epsilon_L(x\wedge y).
\]
Thus \(x\wedge y\in F\), contradicting \(I\cap F=\emptyset\). Therefore
\(I\) is prime. Since every ideal occurring as the first component of a
comaximal pair is prime, Lemma~\ref{ComaximalPrimeCriterion} implies that
\(L\) is distributive.
\end{proof}

For distributive $L$, the sets $d_L(x)$ are compact in $spec(L)$. The next proposition shows that, once the topologies on $\mathfrak{M}(L)$ are defined, the same compactness holds for $\delta_L(x)$ for arbitrary $L$.

\begin{Prop}\label{DelComp}
If $V$ is a non-empty subset of $L$ such that
\[
\delta_L(x)\subseteq \bigcup_{y\in V} \delta_L(y),
\]
then there exists a finite non-empty subset $V_1\subseteq V$ such that
\[
\delta_L(x)\subseteq \bigcup_{y\in V_1} \delta_L(y).
\]
\end{Prop}

\begin{proof}
If $(V]\cap[x)=\emptyset$, then by Proposition~\ref{Com}, there exists a comaximal pair $(I,F)$ such that $V\subseteq I$ and $x\in F$. Hence $(I,F)\in\delta_L(x)\setminus\bigcup_{y\in V}\delta_L(y)$, contradicting the hypothesis. Thus $(V]\cap[x)\neq\emptyset$.  
Hence there exists $z\in (V]\cap[x)$, and thus $v_1,\ldots,v_n\in V$ with
\[
x\le z\le v_1\vee\cdots\vee v_n.
\]
Applying $\delta_L$, a sup-semilattice homomorphism, we obtain
\[
\delta_L(x)\subseteq \delta_L(z)\subseteq \delta_L(v_1)\cup\cdots\cup\delta_L(v_n).
\]
Let $V_1=\{v_1,\ldots,v_n\}$.
\end{proof}

\begin{Lem}
For any $x\in L$, we have $\epsilon_L(x)\subseteq \delta_L(x)$.
\end{Lem}

\begin{proof}
Let \((I,F)\in\epsilon_L(x)\). Then \(x\in F\). Since
\(I\cap F=\emptyset\), we have \(x\notin I\). Hence
\((I,F)\in\delta_L(x)\).
\end{proof}

As a first application, we state one of the most subtle and important properties of this work, connecting the behavior of $\delta_L$ and $\epsilon_L$. We call it the \emph{generalized Balbes--Dwinger property}, since it implies the Balbes--Dwinger property stated in Section~2 when identifying $\delta_L$ (or $\epsilon_L$) with $d_L$ in the distributive case. Later we interpret this as a bitopological property.

\begin{Th}[\textbf{Generalized Balbes--Dwinger property}]\label{GBD}
If $V$ and $W$ are non-empty subsets of $L$ such that
\[
\bigcap_{x\in V} \epsilon_L(x) \subseteq \bigcup_{y\in W} \delta_L(y),
\]
then there exist finite non-empty subsets $V_1\subseteq V$ and $W_1\subseteq W$ such that
\[
\bigcap_{x\in V_1} \epsilon_L(x) \subseteq \bigcup_{y\in W_1} \delta_L(y).
\]
\end{Th}

\begin{proof}
If $(W]\cap[V)=\emptyset$, then by Proposition~\ref{Com}, there exists a comaximal pair $(I,F)$ such that $W\subseteq I$ and $V\subseteq F$. Hence
\[
(I,F)\in \bigcap_{x\in V}\epsilon_L(x)\setminus \bigcup_{y\in W}\delta_L(y),
\]
a contradiction. Thus $(W]\cap[V)\neq\emptyset$.  
Hence there exists $z\in (W]\cap[V)$, so $w_1,\ldots,w_n\in W$ and $v_1,\ldots,v_m\in V$ with
\[
v_1\wedge\cdots\wedge v_m \le z \le w_1\vee\cdots\vee w_n.
\]
Applying $\delta_L$ and $\epsilon_L$ gives
\[
\delta_L(z)\subseteq \delta_L(w_1)\cup\cdots\cup\delta_L(w_n),
\qquad
\epsilon_L(v_1)\cap\cdots\cap\epsilon_L(v_m)\subseteq \epsilon_L(z).
\]
By the previous lemma,
\[
\epsilon_L(v_1)\cap\cdots\cap\epsilon_L(v_m)\subseteq \delta_L(w_1)\cup\cdots\cup\delta_L(w_n).
\]
Let $V_1=\{v_1,\ldots,v_m\}$ and $W_1=\{w_1,\ldots,w_n\}$.
\end{proof}

\subsection{The topologies}
We now have the necessary tools to define two topologies on the set $\mathfrak{M}(L)$. The elements of $Im(\epsilon_L)$ cover $\mathfrak{M}(L)$, and since $\epsilon_L$ is an inf-semilattice homomorphism, $Im(\epsilon_L)$ is a basis, closed under finite intersections, for a topology on $\mathfrak{M}(L)$. We denote this topology by $\sigma_L$ and call it the \emph{inf-topology} on $\mathfrak{M}(L)$. 

In contrast, $Im(\delta_L)$ is in general only a subbasis for a topology on $\mathfrak{M}(L)$. We denote by $\tau_L$ the topology generated by this subbasis, and call it the \emph{sup-topology} on $\mathfrak{M}(L)$.

\begin{Def}
The \emph{bitopological spectrum} of $L$, denoted $spec_B(L)$, is the bitopological space
\[
spec_B(L) = (\mathfrak{M}(L),\tau_L,\sigma_L).
\]
\end{Def}

Given this bitopological space for an arbitrary lattice $L$, the maps $\delta_L$ and $\epsilon_L$ induce lattice embeddings of $L$ into the lattices of $\tau_L$-open and $\sigma_L$-open subsets, respectively. In particular, $\delta_L$ identifies $L$ with a sub-suplattice of $\tau_L$, and $\epsilon_L$ identifies $L$ with a sub-inflattice of $\sigma_L$. Our aim is to recover (at least) one of these copies of $L$. To do so, we study certain bitopological properties of $spec_B(L)$ and develop techniques adapted to this setting.

\begin{Prop}
For any $x\in L$, the set $\delta_L(x)$ is $\tau_L$-compact.
\end{Prop}

\begin{proof}
This follows from Proposition~\ref{DelComp} together with the Alexander subbasis lemma.
\end{proof}

\begin{Def}[\cite{Kelly}]
A bitopological space $(X,\tau,\sigma)$ is \emph{pairwise $T_0$} if for any distinct points $x,y\in X$ there exists either $U\in\tau$ with $x\in U$ and $y\notin U$, or $V\in\sigma$ with $x\notin V$ and $y\in V$.
\end{Def}

We will show that $spec_B(L)$ is pairwise $T_0$. To this end, we characterize the topological preorders induced by $\tau_L$ and $\sigma_L$, and relate the pairwise $T_0$ property to a property of preordered sets.

\begin{Def}[\cite{Urquhart}]
A set $X$ equipped with two preorders $\leq_1$ and $\leq_2$ is \emph{pairwise ordered} if, for all $x,y\in X$,
\[
x\leq_1 y \text{ and } y\leq_2 x \quad\Rightarrow\quad x=y.
\]
\end{Def}

Let $(X,\tau,\sigma)$ be any bitopological space, and denote by $\leq_\tau$ and $\leq_\sigma$ the topological preorders induced by $\tau$ and $\sigma$, respectively (i.e.\ $x\leq_\tau y$ iff $x\in\overline{\{y\}}^{\tau}$, and similarly for $\sigma$). From the definitions, the following are equivalent:
\begin{enumerate}
    \item[(i)] $(X,\tau,\sigma)$ is pairwise $T_0$,
    \item[(ii)] $(X,\leq_\tau,\leq_\sigma)$ is a pairwise ordered set.
\end{enumerate}
Thus, to prove that $spec_B(L)$ is pairwise $T_0$, it suffices to show that it is pairwise ordered with respect to the topological preorders induced by $\tau_L$ and $\sigma_L$. We first require a simple lemma.

\begin{Lem}\label{RedSub}
Let $(X,\tau)$ be a topological space with subbasis $S$, and let $x,y\in X$. If for every $U\in S$,
\[
x\in U \ \Rightarrow\ y\in U,
\]
then $x\in\overline{\{y\}}^{\tau}$, i.e.\ $x\leq_\tau y$.
\end{Lem}

Denote by $\leq_{\tau}$ the topological preorder on $spec_B(L)$ induced by $\tau_L$, and by $\leq_{\sigma}$ the one induced by $\sigma_L$.

\begin{Prop}\label{Ord}
Let $(I,F)$ and $(J,G)$ be comaximal pairs of $L$. Then:
\begin{enumerate}
    \item[(i)] $(I,F)\leq_{\tau} (J,G)$ if and only if $J\subseteq I$.
    \item[(ii)] $(I,F)\leq_{\sigma} (J,G)$ if and only if $F\subseteq G$.
\end{enumerate}
\end{Prop}

\begin{proof}
\begin{enumerate}
    \item[(i)] Suppose $(I,F)\leq_{\tau} (J,G)$ but $J\not\subseteq I$. Then there exists $x\in J\setminus I$. Since $x\notin I$, we have $(I,F)\in\delta_L(x)$; but $x\in J$ implies $(J,G)\notin\delta_L(x)$, contradicting $(I,F)\leq_{\tau} (J,G)$. Hence $J\subseteq I$.

    Conversely, assume $J\subseteq I$. For any $x\in L$, if $(I,F)\in\delta_L(x)$, then $x\notin I$. Since $J\subseteq I$, this implies $x\notin J$, and hence $(J,G)\in\delta_L(x)$. Thus every subbasic $\tau_L$-open set containing $(I,F)$ also contains $(J,G)$, and Lemma~\ref{RedSub} yields $(I,F)\leq_{\tau} (J,G)$.

    \item[(ii)] The proof is analogous. If $(I,F)\leq_{\sigma} (J,G)$ but $F\not\subseteq G$, choose $x\in F\setminus G$. Then $(I,F)\in\epsilon_L(x)$ and $(J,G)\notin\epsilon_L(x)$, a contradiction. Conversely, if $F\subseteq G$, any subbasic $\sigma_L$-open set $\epsilon_L(x)$ containing $(I,F)$ (i.e.\ with $x\in F$) also contains $(J,G)$, since $x\in F\subseteq G$. Again Lemma~\ref{RedSub} applies.
\end{enumerate}
\end{proof}

\begin{Cor}
\(spec_B(L)\) is a pairwise \(T_0\) bitopological space.
\end{Cor}

\begin{proof}
Let
\[
    p=(I,F), \qquad q=(J,G)
\]
be two distinct comaximal pairs of \(L\). Suppose, toward a contradiction,
that there is no \(\tau_L\)-open set containing \(p\) and not containing
\(q\), and no \(\sigma_L\)-open set containing \(q\) and not containing
\(p\). Then
\[
    p\leq_{\tau_L}q
    \qquad\text{and}\qquad
    q\leq_{\sigma_L}p.
\]
By Proposition~\ref{Ord}, this means
\[
    J\subseteq I
    \qquad\text{and}\qquad
    G\subseteq F.
\]

If \(J\subsetneq I\), then \(I\) is an ideal properly containing \(J\). Since
\((J,G)\) is a comaximal pair, we must have
\[
    I\cap G\neq\emptyset.
\]
But \(G\subseteq F\), and therefore
\[
    I\cap F\neq\emptyset,
\]
contradicting the fact that \((I,F)\) is a comaximal pair. Hence \(J=I\).

Similarly, if \(G\subsetneq F\), then \(F\) is a filter properly containing
\(G\). Since \((J,G)\) is a comaximal pair, we must have
\[
    J\cap F\neq\emptyset.
\]
Since \(J=I\), this contradicts \(I\cap F=\emptyset\). Hence \(G=F\).

Therefore \((I,F)=(J,G)\), contradicting our assumption that \(p\neq q\).
Thus the required pairwise \(T_0\) separation holds.
\end{proof}

In a certain sense, the sup-topology on $spec_B(L)$ encodes the supremum-semilattice structure of $L$, while the inf-topology encodes the infimum-semilattice structure. The following propositions illustrate how lattice-theoretic properties are reflected in (bi)topological terms.

\begin{Prop}\label{1IffComp}
The following are equivalent:
\begin{enumerate}
    \item[(i)] $L$ has a top element $1$,
    \item[(ii)] $spec_B(L)$ is $\tau_L$-compact.
\end{enumerate}
\end{Prop}

\begin{proof}
Assume $L$ has $1$. Then $spec_B(L)=\delta_L(1)$, which is $\tau_L$-compact by Proposition~\ref{DelComp}.

Conversely, suppose $spec_B(L)$ is $\tau_L$-compact. Since the family $\{\delta_L(x)\mid x\in L\}$ covers $spec_B(L)$, we have
\[
spec_B(L)=\bigcup_{x\in L}\delta_L(x).
\]
By compactness, there exist $x_1,\dots,x_n\in L$ such that
\[
spec_B(L)=\delta_L(x_1)\cup\cdots\cup\delta_L(x_n)
=\delta_L(x_1\vee\cdots\vee x_n).
\]
We claim that $u:=x_1\vee\cdots\vee x_n$ is the maximum of $L$. Suppose not. Then there exists $y\in L$ with $u\not\ge y$, so $(u]\cap [y)=\emptyset$. By Proposition~\ref{Com}, there exists a comaximal pair $(I,F)$ with $u\in I$ and $y\in F$. Thus $(I,F)\notin\delta_L(u)$, contradicting $spec_B(L)=\delta_L(u)$. Hence $u$ is the top element of $L$.
\end{proof}

\begin{Prop}\label{0IffFun}
The following are equivalent:
\begin{enumerate}
    \item[(i)] $L$ has a bottom element $0$,
    \item[(ii)] $\emptyset$ is $\sigma_L$-fundamental.
\end{enumerate}
\end{Prop}

\begin{proof}
First observe that
\[
\bigcap_{x\in L}\epsilon_L(x)=\emptyset,
\]
since no comaximal pair $(I,F)$ can contain every element of $L$ in $F$.

Assume $L$ has $0$. To show that $\emptyset$ is $\sigma_L$-fundamental, it suffices to prove that if $V\subseteq L$ satisfies
\[
\bigcap_{x\in V}\epsilon_L(x)=\emptyset,
\]
then there exists a finite subset $V_1\subseteq V$ with $\bigcap_{x\in V_1}\epsilon_L(x)=\emptyset$. But
\[
\bigcap_{x\in V}\epsilon_L(x)\subseteq \emptyset =\delta_L(0),
\]
and so Theorem~\ref{GBD} applies.

Conversely, suppose $\emptyset$ is $\sigma_L$-fundamental. Since $\bigcap_{x\in L}\epsilon_L(x)=\emptyset$, there exist $x_1,\dots,x_n\in L$ such that
\[
\epsilon_L(x_1)\cap\cdots\cap\epsilon_L(x_n)=\emptyset,
\]
i.e.\ $\epsilon_L(x_1\wedge\cdots\wedge x_n)=\emptyset$. If $u:=x_1\wedge\cdots\wedge x_n$ were not the minimum, there would exist some $y\in L$ with $[y)\cap(u]=\emptyset$, and by Proposition~\ref{Com} we could construct a comaximal pair $(I,F)$ containing $y$ in $F$ but not $u$, contradicting $\epsilon_L(u)=\emptyset$. Hence $u$ is the minimum element of $L$.
\end{proof}

We conclude this section with a bitopological criterion detecting when the
first component of a comaximal pair is a prime ideal. This criterion will be
used to characterize distributivity in terms of the coincidence of the two
point-closures throughout the spectrum.

\begin{Prop}[\textbf{Bitopological criterion for prime ideals}]\label{BCPrimal}
Let \((I,F)\) be a comaximal pair of \(L\). If
\[
\overline{\{(I,F)\}}^{\tau_L}
=
\overline{\{(I,F)\}}^{\sigma_L},
\]
then \(I\) is a prime ideal of \(L\).
\end{Prop}

\begin{proof}
We first show that \(F=I^c\). Since \(I\cap F=\emptyset\), we have
\[
    F\subseteq I^c.
\]

Conversely, let \(x\in I^c\). Then \(x\notin I\), and therefore the ideal
\(I\) and the principal filter \([x)\) are disjoint. By
Proposition~\ref{Com}, there exists a comaximal pair \((J,G)\) such that
\[
    I\subseteq J
    \qquad\text{and}\qquad
    [x)\subseteq G.
\]
In particular, \(x\in G\).

Since \(I\subseteq J\), Proposition~\ref{Ord}(i) implies
\[
    (J,G)\in\overline{\{(I,F)\}}^{\tau_L}.
\]
By the hypothesis on the two closures, we also have
\[
    (J,G)\in\overline{\{(I,F)\}}^{\sigma_L}.
\]
Using Proposition~\ref{Ord}(ii), this means
\[
    G\subseteq F.
\]
Since \(x\in G\), it follows that \(x\in F\). Hence \(I^c\subseteq F\), and
therefore
\[
    F=I^c.
\]

Since \(F\) is a filter and \(F=I^c\), the complement of \(I\) is a filter.
Thus \(I\) is a prime ideal.
\end{proof}

\begin{Def}
Let $(X,\tau,\sigma)$ be a bitopological space. A point $x\in X$ is called \textbf{prime} if
\[
\overline{\{x\}}^{\tau}=\overline{\{x\}}^{\sigma}.
\]
\end{Def}

\begin{Th}\label{PrimePointsDist}
A lattice \(L\) is distributive if and only if every point of \(spec_B(L)\)
is prime.
\end{Th}

\begin{proof}
Suppose first that \(L\) is distributive. By Proposition~\ref{TopDis}, we
have
\[
    \delta_L=\epsilon_L.
\]
Hence
\[
    \operatorname{Im}(\delta_L)=\operatorname{Im}(\epsilon_L).
\]
Since \(\operatorname{Im}(\epsilon_L)\) is a basis for \(\sigma_L\), this
common family is closed under finite intersections. Therefore the topology
\(\tau_L\), generated by \(\operatorname{Im}(\delta_L)\) as a subbasis, and
the topology \(\sigma_L\), generated by \(\operatorname{Im}(\epsilon_L)\) as
a basis, coincide. Thus
\[
    \tau_L=\sigma_L.
\]
Consequently, the two closures of every point of \(spec_B(L)\) coincide, so
every point of \(spec_B(L)\) is prime.

Conversely, suppose that every point of \(spec_B(L)\) is prime. Let
\((I,F)\) be a comaximal pair of \(L\). Since \((I,F)\) is a prime point, we
have
\[
\overline{\{(I,F)\}}^{\tau_L}
=
\overline{\{(I,F)\}}^{\sigma_L}.
\]
By Proposition~\ref{BCPrimal}, \(I\) is a prime ideal. Hence every ideal
which occurs as the first component of a comaximal pair of \(L\) is prime.
By Lemma~\ref{ComaximalPrimeCriterion}, \(L\) is distributive.
\end{proof}

\subsection{Reconstruction theory}
In this section we present the way to recover an arbitrary lattice $L$ from its bitopological spectrum. To this end, we extend the techniques developed by Balbes and Dwinger, adapting them to our setting, and we also adjust part of Urquhart's reconstruction method \cite{Urquhart} for recovering a lattice from its Urquhart spectrum. We know that there is a copy of $L$ inside the sup-topology $\tau_L$, and an isomorphic copy of $L$ inside the inf-topology $\sigma_L$, of its spectrum. In a sense, the ``upper'' part of the lattice structure is encoded in $\tau_L$, while the ``lower'' part is encoded in $\sigma_L$. Using certain bitopological properties and the transition maps induced by the
two topological preorders, we will ``glue'' these two copies.

Let $(X,\leq_1,\leq_2)$ be a doubly ordered set (i.e.\ a set with two order relations). For $i\in\{1,2\}$, define
\begin{align*}
    \downarrow_i(A) &= \{x\in X \mid \text{there exists } a\in A \text{ such that } x\leq_i a\},\\[2mm]
    \uparrow_i(A)  &= \{x\in X \mid \text{there exists } a\in A \text{ such that } a\leq_i x\}.
\end{align*}
and let $C_i(X)$ denote the collection of all $\leq_i$-increasing subsets of $X$. Define
\begin{eqnarray*}
    d : C_1(X) &\longrightarrow& C_2(X)\\
    A &\longmapsto& \bigl(\downarrow_2(A^c)\bigr)^c,
\end{eqnarray*}
and
\begin{eqnarray*}
    i : C_2(X) &\longrightarrow& C_1(X)\\
    A &\longmapsto& \uparrow_1(A).
\end{eqnarray*}
The maps \(i\) and \(d\) are order-preserving. Moreover, \(d\) preserves
arbitrary intersections and \(i\) preserves arbitrary unions.

Since any bitopological space $(X,\tau,\sigma)$ induces a doubly preordered set $(X,\leq_\tau,\leq_\sigma)$, it makes sense to define $i$ and $d$ in the bitopological context.

For a bitopological space $(X,\tau,\sigma)$, we denote by $i_X$ and $d_X$ the
transition maps introduced above, or simply by \(i\) and \(d\) when the
context is clear. We refer to \(i\) and \(d\) as the \textbf{transition
functions} of \(X\).
\begin{Def}
Let $(X,\tau,\sigma)$ be a bitopological space.
\begin{enumerate}
    \item[(i)] A $\leq_\tau$-increasing subset $A\subseteq X$ is \emph{stable} if $i\circ d(A)=A$.
    \item[(ii)] A $\leq_\sigma$-increasing subset $B\subseteq X$ is \emph{co-stable} if $d\circ i(B)=B$.
\end{enumerate}
\end{Def}

Note that any open set is increasing with respect to its corresponding topological preorder.

\medskip

We return now to lattices. Let $L$ be a lattice. The idea behind introducing $i$ and $d$ is to relate the functions $\delta_L$ and $\epsilon_L$.

\begin{Prop}
Let $x\in L$. Then:
\begin{enumerate}
    \item[(i)] $d(\delta_L(x))=\epsilon_L(x)$,
    \item[(ii)] $i(\epsilon_L(x))=\delta_L(x)$,
    \item[(iii)] $\delta_L(x)$ is a stable subset of $spec_B(L)$,
    \item[(iv)] $\epsilon_L(x)$ is a co-stable subset of $spec_B(L)$.
\end{enumerate}
\end{Prop}

\begin{proof}
\begin{enumerate}
    \item[(i)] Let $(I,F)\in d(\delta_L(x))$. Suppose $(I,F)\notin\epsilon_L(x)$, i.e.\ $x\notin F$. By Proposition~\ref{Com}, there exists a comaximal pair $(J,G)$ of $L$ such that $x\in J$ and $F\subseteq G$. Then $(I,F)\leq_\sigma (J,G)$, and since $d(\delta_L(x))$ is $\leq_\sigma$-increasing, we have $(J,G)\in d(\delta_L(x))$. On the other hand, \(x\in J\) implies \((J,G)\in\delta_L(x)^c\). Since
    \((I,F)\leq_\sigma (J,G)\), it follows that
    \[
        (I,F)\in\downarrow_\sigma(\delta_L(x)^c),
    \]
    contradicting the fact that
    \[
        (I,F)\in d(\delta_L(x))
        =
        \bigl(\downarrow_\sigma(\delta_L(x)^c)\bigr)^c.
    \]
    
    Conversely, let $(I,F)\in\epsilon_L(x)$, so $x\in F$. Suppose $(J,G)\in\delta_L(x)^c$ and $F\subseteq G$. Then $x\in G$, contradicting $(J,G)\in\delta_L(x)^c$ (which means $x\in J$). Thus for any $(J,G)\in\delta_L(x)^c$, we cannot have $(I,F)\leq_\sigma (J,G)$, and by definition this means $(I,F)\in d(\delta_L(x))$.

    \item[(ii)] Let $(I,F)\in i(\epsilon_L(x))$. Then there is $(J,G)\in\epsilon_L(x)$ with $I\subseteq J$. Since $x\in G$ and $(J,G)$ is comaximal, we have $x\notin J$, so $x\notin I$, and hence $(I,F)\in\delta_L(x)$.

    Conversely, let $(I,F)\in\delta_L(x)$, so $x\notin I$. By Proposition~\ref{Com}, there exists a comaximal pair $(J,G)$ with $I\subseteq J$ and $x\in G$. Then $(J,G)\in\epsilon_L(x)$ and, by definition of $i$, $(I,F)\in i(\epsilon_L(x))$.

    \item[(iii)] Since $d(\delta_L(x))=\epsilon_L(x)$ and $i(\epsilon_L(x))=\delta_L(x)$, we have
    \[
      i\circ d(\delta_L(x)) = i(\epsilon_L(x)) = \delta_L(x),
    \]
    so $\delta_L(x)$ is stable.

    \item[(iv)] Similarly, $d\circ i(\epsilon_L(x)) = d(\delta_L(x)) = \epsilon_L(x)$, so $\epsilon_L(x)$ is co-stable.
\end{enumerate}
\end{proof}

\begin{Def}
Let $(X,\tau,\sigma)$ be a bitopological space. A non-empty subset $A\subseteq X$ is \textbf{essential} if
\begin{enumerate}
    \item $A$ is $\tau$-compact,
    \item $d(A)$ is $\sigma$-open,
    \item $A$ is stable.
\end{enumerate}
Let \(\mathfrak{E}^{+}(X)\) denote the collection of non-empty essential
subsets of \(X\). The empty set is said to be \textbf{\(\sigma\)-fundamental}
if every family
\[
    \{d(A)\mid A\in\mathcal{V}\},
    \qquad \mathcal{V}\subseteq \mathfrak{E}^{+}(X),
\]
with the finite intersection property has non-empty intersection.

The empty set is \textbf{essential} if it is \(\sigma\)-fundamental.
\end{Def}

Essential subsets generalize fundamental subsets to the bitopological setting.
Indeed, when the two topologies coincide, the transition functions restrict
to the identity on the open subsets of the common topology. In that case, the
collections of fundamental and essential subsets coincide.

As mentioned in Section~2, Balbes, Dwinger and Stone recover any distributive lattice from its (Stone) spectrum by considering the fundamental subsets of that topological space. In our setting, an arbitrary lattice will be recovered from the essential subsets of its bitopological spectrum.

\begin{Th}\label{RepTh}
Every essential subset of $spec_B(L)$ is of the form $\delta_L(x)$ for some $x\in L$.
\end{Th}

\begin{proof}
Let $A$ be a non-empty essential subset of $spec_B(L)$. Since $d(A)$ is $\sigma_L$-open, there is a subset $W\subseteq L$ such that
\[
d(A) = \bigcup_{x\in W}\epsilon_L(x).
\]
Then
\begin{eqnarray*}
A &=& i(d(A)) \qquad\text{(since $A$ is stable)}\\
  &=& i\Bigl(\bigcup_{x\in W}\epsilon_L(x)\Bigr)\\
  &=& \bigcup_{x\in W} i(\epsilon_L(x)) \qquad\text{($i$ preserves arbitrary unions)}\\
  &=& \bigcup_{x\in W} \delta_L(x).
\end{eqnarray*}
Since $A$ is $\tau_L$-compact and each $\delta_L(x)$ is $\tau_L$-open, there exist $w_1,\dots,w_n\in W$ such that
\[
A = \delta_L(w_1)\cup\cdots\cup\delta_L(w_n)
  = \delta_L(w_1\vee\cdots\vee w_n),
\]
using that $\delta_L$ is a sup-semilattice homomorphism.

If the empty set is essential, then by definition it is $\sigma_L$-fundamental, and Proposition~\ref{0IffFun} implies that $L$ has $0$. Hence $\emptyset = \delta_L(0)$.
\end{proof}

\section{An extension of the Stone duality theorem}
\subsection{Functorial behaviour of the construction}
In the case of the construction of the (Stone) spectrum of a distributive lattice, not every homomorphism between distributive lattices preserves prime ideals by inverse preimage. The class of homomorphisms with this property are called \textit{proper}. In our construction, we have a similar failure of having full functoriality: it is not true that every lattice homomorphism preserves maximal copairs by inverse image.

\begin{Def}
A lattice homomorphism $f:L\rightarrow M$ is \textbf{quasi-proper} if $(f^{-1}(I),f^{-1}(F))$ is a comaximal pair of $L$ for every comaximal pair $(I,F)$ of $M$.
\end{Def}

Any quasi-proper homomorphism is proper, and for distributive lattices the notions of quasi-proper and proper coincide. However, in the non-distributive case there exist proper homomorphisms that are not quasi-proper. (For example, the inclusion of any two-element chain into $\mathfrak{M}_5$ is vacuously proper, but not all of these inclusions are quasi-proper.)

\begin{Def}
A bicontinuous function $f:(X,\tau_1,\sigma_1)\rightarrow(Y,\tau_2,\sigma_2)$ is \textbf{strongly bicontinuous} if it preserves essential sets under inverse image.
\end{Def}

In this work we will not be able to represent all lattice homomorphisms, but we do represent all objects in the category of lattices. Let $\mathfrak{R}_{qp}$ denote the subcategory of all lattices with quasi-proper homomorphisms. This is the category that we will dualize to a suitable category of bitopological spaces.

\begin{Prop}\label{Func}
If $f:L\rightarrow N$ is a quasi-proper lattice homomorphism, then the map
\[
\begin{array}{rcl}
spec_B(f):spec_B(N) &\longrightarrow& spec_B(L)\\[1mm]
(I,F) &\longmapsto& (f^{-1}(I),f^{-1}(F)),
\end{array}
\]
is strongly bicontinuous.
\end{Prop}

\begin{proof}
Let \(x\in L\). Then:
\begin{align*}
(I,F)\in spec_B(f)^{-1}(\delta_L(x))
&\Longleftrightarrow spec_B(f)(I,F)\in \delta_L(x)\\
&\Longleftrightarrow x\notin f^{-1}(I)\\
&\Longleftrightarrow f(x)\notin I\\
&\Longleftrightarrow (I,F)\in \delta_N(f(x)).
\end{align*}
Thus
\[
    spec_B(f)^{-1}(\delta_L(x))=\delta_N(f(x)).
\]
Similarly,
\[
    spec_B(f)^{-1}(\epsilon_L(x))=\epsilon_N(f(x)).
\]

Since the sets \(\delta_L(x)\) form a subbasis for \(\tau_L\), the first
identity shows that \(spec_B(f)\) is continuous with respect to the
sup-topologies. Since the sets \(\epsilon_L(x)\) form a basis for
\(\sigma_L\), the second identity shows that \(spec_B(f)\) is continuous
with respect to the inf-topologies. Hence \(spec_B(f)\) is bicontinuous.

It remains to show that inverse images of essential subsets are essential.
Let \(A\subseteq spec_B(L)\) be essential. By Theorem~\ref{RepTh}, there
exists \(x\in L\) such that
\[
    A=\delta_L(x).
\]
Therefore
\[
    spec_B(f)^{-1}(A)
    =
    spec_B(f)^{-1}(\delta_L(x))
    =
    \delta_N(f(x)).
\]
By the results preceding Theorem~\ref{RepTh}, the set
\(\delta_N(f(x))\) is essential in \(spec_B(N)\). Hence \(spec_B(f)\) preserves essential
subsets under inverse image.

Consequently, \(spec_B(f)\) is strongly bicontinuous.
\end{proof}
Therefore, $spec_B$ may be regarded as a functor from the category $\mathfrak{R}_{qp}$ to the category of bitopological spaces with strongly bicontinuous maps.

\subsection{Pairwise Balbes--Dwinger spaces}
In this section we will characterize the essential image of the functor $spec_B$, i.e., we will find conditions for a bitopological space $(X,\tau,\sigma)$ which guarantees that $X$ is bihomeomorphic to the bitopological spectrum of some lattice. These conditions extend the properties of a Balbes--Dwinger space.

\begin{Def}
Let $X$ be a set and let $f:\wp(X)\rightarrow\wp(X)$ be a function. A collection
$\mathcal{H}$ of subsets of \(X\) is \textbf{\(f\)-birreducible} if for any
families \(V,W\subseteq\mathcal{H}\) such that
\[
\bigcap_{A\in V} f(A)\subseteq \bigcup_{B\in W}B,
\]
there exist finite subsets \(V_0\subseteq V\) and \(W_0\subseteq W\), with
\(V_0\) non-empty whenever \(V\) is non-empty, such that
\[
\bigcap_{A\in V_0} f(A)\subseteq \bigcup_{B\in W_0}B.
\]
As usual, an empty union is interpreted as \(\emptyset\), and an empty
intersection is interpreted as \(X\).
\end{Def}

\begin{Def}
A bitopological space $(X,\tau,\sigma)$ is a \textbf{pairwise Balbes-Dwinger space} if:
\begin{enumerate}
    \item[(i)] $(X,\tau,\sigma)$ is pairwise $T_0$.
    \item[(ii)] The collection of essential sets is a sub-basis for $\tau$.
    \item[(iii)] The set $\{d(A)\mid A \text{ essential}\}$ is a basis for $\sigma$, closed under finite intersections.
    \item[(iv)] If $A$ and $B$ are essential, then $A\cup B$ and $i(d(A\cap B))$ are essential.
    \item[(v)] The collection of essential subsets of $X$ is $d$-birreducible.
\end{enumerate}
For a pairwise Balbes-Dwinger space $(X,\tau,\sigma)$ we denote by $\mathfrak{E}(X)$ the collection of essential subsets of $X$.
\end{Def}

The preceding results show that, for every lattice \(L\), the bitopological
space \(spec_B(L)\) is a pairwise Balbes--Dwinger space.

\begin{Prop}
If $(X,\tau,\sigma)$ is a pairwise Balbes-Dwinger space, then $(\mathfrak{E}(X),\subseteq)$ is a lattice with operations:
\begin{enumerate}
    \item[(i)] $A\vee B = A\cup B$,
    \item[(ii)] $A\wedge B = i(d(A\cap B))$.
\end{enumerate}
\end{Prop}

\begin{proof}
By axiom \((iv)\) in the definition of pairwise Balbes--Dwinger space,
\(A\cup B\) is essential. It is clearly the least upper bound of \(A\) and
\(B\), so
\[
    A\vee B=A\cup B.
\]

Now set
\[
    M=i(d(A\cap B)).
\]
Again by axiom \((iv)\), \(M\) is essential. We show that \(M\) is the
greatest lower bound of \(A\) and \(B\).

Since \(A\cap B\subseteq A\), monotonicity of \(d\) and \(i\) gives
\[
    i(d(A\cap B))\subseteq i(d(A)).
\]
Because \(A\) is stable, \(i(d(A))=A\). Hence \(M\subseteq A\). Similarly,
\(M\subseteq B\). Thus \(M\) is a lower bound of \(A\) and \(B\).

Let \(C\in\mathfrak{E}(X)\) be such that \(C\subseteq A\) and
\(C\subseteq B\). Then \(C\subseteq A\cap B\). By monotonicity,
\[
    d(C)\subseteq d(A\cap B),
\]
and hence
\[
    i(d(C))\subseteq i(d(A\cap B)).
\]
Since \(C\) is stable, \(i(d(C))=C\). Therefore
\[
    C\subseteq i(d(A\cap B))=M.
\]
Thus \(M\) is the greatest lower bound of \(A\) and \(B\). Consequently,
\[
    A\wedge B=i(d(A\cap B)).
\]
\end{proof}

Now we relate $X$ with $spec_B(\mathfrak{E}(X))$.

\begin{Lem}\label{dEmptyIntersection}
If \((X,\tau,\sigma)\) is a pairwise Balbes--Dwinger space, then
\[
\bigcap_{A\in\mathfrak{E}(X)} d(A)=\emptyset.
\]
\end{Lem}

\begin{proof}
Suppose, toward a contradiction, that there exists \(x\in X\) such that
\[
    x\in \bigcap_{A\in\mathfrak{E}(X)} d(A).
\]
Then every family of basic \(\sigma\)-open sets of the form \(d(A)\), with
\(A\in\mathfrak{E}(X)\), has non-empty intersection, since all of them
contain \(x\). Because the family
\[
    \{d(A)\mid A\in\mathfrak{E}(X)\}
\]
is a basis for \(\sigma\), closed under finite intersections, this means that
\(\emptyset\) is \(\sigma\)-fundamental. Hence \(\emptyset\) is essential.

But then \(\emptyset\in\mathfrak{E}(X)\), and so
\[
    x\in d(\emptyset).
\]
This is impossible, since
\[
    d(\emptyset)
    =
    \left(\downarrow_\sigma(X)\right)^c
    =
    \emptyset.
\]
Therefore
\[
    \bigcap_{A\in\mathfrak{E}(X)} d(A)=\emptyset.
\]
\end{proof}

\begin{Lem}\label{EssentialEmptyIntersection}
If \((X,\tau,\sigma)\) is a pairwise Balbes--Dwinger space, then
\[
\bigcap_{A\in\mathfrak{E}(X)} A=\emptyset.
\]
\end{Lem}

\begin{proof}
Suppose, toward a contradiction, that there exists \(x\in X\) such that
\[
    x\in \bigcap_{A\in\mathfrak{E}(X)} A.
\]
We prove that \(\emptyset\) is \(\sigma\)-fundamental.

Since the sets \(d(A)\), with \(A\in\mathfrak{E}(X)\), form a basis for
\(\sigma\) closed under finite intersections, it is enough to verify the
finite intersection property for families of such basic open sets.

Let \(\mathcal{V}\subseteq\mathfrak{E}(X)\) be such that the family
\[
    \{d(A)\mid A\in\mathcal{V}\}
\]
has the finite intersection property. Suppose, toward a contradiction, that
\[
    \bigcap_{A\in\mathcal{V}} d(A)=\emptyset.
\]
Equivalently,
\[
    \bigcap_{A\in\mathcal{V}} d(A)
    \subseteq
    \bigcup_{B\in\emptyset}B.
\]
By \(d\)-birreducibility, applied to the family \(\mathcal{V}\) and to the
empty family on the right-hand side, there exists a finite subset
\(\mathcal{V}_0\subseteq\mathcal{V}\) such that
\[
    \bigcap_{A\in\mathcal{V}_0} d(A)=\emptyset.
\]
This contradicts the finite intersection property. Hence
\[
    \bigcap_{A\in\mathcal{V}} d(A)\neq\emptyset.
\]
Therefore \(\emptyset\) is \(\sigma\)-fundamental, and hence \(\emptyset\) is
essential.

But this contradicts the assumption that \(x\) belongs to every essential
subset, since \(x\notin\emptyset\). Therefore
\[
    \bigcap_{A\in\mathfrak{E}(X)} A=\emptyset.
\]
\end{proof}

\begin{Lem}\label{dMeet}
Let \(X\) be a pairwise Balbes--Dwinger space. If \(A,B\in\mathfrak{E}(X)\), then
\[
d(A\wedge B)=d(A)\cap d(B).
\]
\end{Lem}

\begin{proof}
Since the family \(\{d(C)\mid C\in\mathfrak{E}(X)\}\) is closed under finite
intersections, there exists \(C\in\mathfrak{E}(X)\) such that
\[
d(C)=d(A)\cap d(B).
\]
Since \(C\) is stable, \(C=i(d(C))\). Hence
\[
C=i(d(C))\subseteq i(d(A))=A
\qquad\text{and}\qquad
C=i(d(C))\subseteq i(d(B))=B,
\]
so \(C\subseteq A\wedge B\). Therefore
\[
d(C)\subseteq d(A\wedge B).
\]
On the other hand, since \(A\wedge B\subseteq A\) and \(A\wedge B\subseteq B\), we have
\[
d(A\wedge B)\subseteq d(A)\cap d(B)=d(C).
\]
Thus \(d(A\wedge B)=d(A)\cap d(B)\).
\end{proof}

\begin{Lem}\label{CharMaxPair}
Let $(X,\tau,\sigma)$ be a pairwise Balbes-Dwinger space. For any $x\in X$ define:
\begin{enumerate}
    \item[(i)] $I(x)=\{A\in\mathfrak{E}(X)\mid x\notin A\}$,
    \item[(ii)] $F(x)=\{A\in\mathfrak{E}(X)\mid x\in d(A)\}$.
\end{enumerate}
Then all comaximal pairs of $\mathfrak{E}(X)$ are of the form $(I(x),F(x))$ for some $x$, and $(I(x),F(x))\neq(I(y),F(y))$ whenever $x\neq y$.
\end{Lem}

\begin{proof}
We divide the proof into several steps.

\begin{enumerate}
    \item[(i)] First, \(I(x)\) is an ideal of \(\mathfrak{E}(X)\). Indeed,
    suppose that \(A\subseteq B\), with \(B\in I(x)\) and
    \(A\in\mathfrak{E}(X)\). Since \(x\notin B\), we also have
    \(x\notin A\), and hence \(A\in I(x)\). Moreover, if
    \(A,B\in I(x)\), then \(x\notin A\cup B\). Since
    \(A\vee B=A\cup B\), it follows that \(A\vee B\in I(x)\). Finally,
    \(I(x)\neq\emptyset\), because
    \[
        \bigcap_{A\in\mathfrak{E}(X)}A=\emptyset
    \]
    by Lemma~\ref{EssentialEmptyIntersection}.

    \item[(ii)] Similarly, \(F(x)\) is a filter of \(\mathfrak{E}(X)\).
    Suppose first that \(A\in F(x)\), \(A\subseteq B\), and
    \(B\in\mathfrak{E}(X)\). Since \(d\) is order-preserving, we have
    \(d(A)\subseteq d(B)\). Thus \(x\in d(B)\), and therefore
    \(B\in F(x)\).

    Now let \(A,B\in F(x)\). Then \(x\in d(A)\cap d(B)\). By
    Lemma~\ref{dMeet},
    \[
        d(A)\cap d(B)=d(A\wedge B).
    \]
    Hence \(x\in d(A\wedge B)\), and so \(A\wedge B\in F(x)\).

    Finally, \(F(x)\neq\emptyset\), since the family
    \[
        \{d(A)\mid A\in\mathfrak{E}(X)\}
    \]
    is a basis for \(\sigma\). Hence, for the point \(x\), there exists
    \(A\in\mathfrak{E}(X)\) such that \(x\in d(A)\), and therefore
    \(A\in F(x)\).

    \item[(iii)] We have \(I(x)\cap F(x)=\emptyset\). Indeed, if
    \(A\in I(x)\cap F(x)\), then \(x\notin A\) and \(x\in d(A)\). This is impossible, because \(A\) is stable and therefore
\[
    d(A)\subseteq i(d(A))=A.
\]
Hence
    \(I(x)\cap F(x)=\emptyset\).

    \item[(iv)] We now show that every comaximal pair of
    \(\mathfrak{E}(X)\) is of the form \((I(x),F(x))\). Let \((I,F)\) be
    a comaximal pair of \(\mathfrak{E}(X)\). We claim that
    \[
        \bigcap_{A\in F}d(A)\not\subseteq \bigcup_{B\in I}B.
    \]
    Suppose, toward a contradiction, that
    \[
        \bigcap_{A\in F}d(A)\subseteq \bigcup_{B\in I}B.
    \]
    Since the collection of essential subsets of \(X\) is
    \(d\)-birreducible, there exist
    \(A_1,\ldots,A_n\in F\) and \(B_1,\ldots,B_m\in I\) such that
    \[
        d(A_1)\cap\cdots\cap d(A_n)
        \subseteq
        B_1\cup\cdots\cup B_m.
    \]
    Since \(d\) preserves arbitrary intersections, this gives
    \[
        d(A_1\cap\cdots\cap A_n)
        \subseteq
        B_1\cup\cdots\cup B_m.
    \]
    The set \(B_1\cup\cdots\cup B_m\) is \(\tau\)-open, and therefore it
    is \(\leq_{\tau}\)-increasing. Hence
    \[
        i\bigl(d(A_1\cap\cdots\cap A_n)\bigr)
        \subseteq
        B_1\cup\cdots\cup B_m.
    \]
    The left-hand side is the finite meet
    \(A_1\wedge\cdots\wedge A_n\), while the right-hand side is the finite
    join \(B_1\vee\cdots\vee B_m\). Thus
    \[
        A_1\wedge\cdots\wedge A_n
        \subseteq
        B_1\vee\cdots\vee B_m.
    \]
    Since \(F\) is a filter, \(A_1\wedge\cdots\wedge A_n\in F\). Since
    \(I\) is an ideal, \(B_1\vee\cdots\vee B_m\in I\). The preceding
    inclusion then implies
    \[
        A_1\wedge\cdots\wedge A_n\in I,
    \]
    contradicting \(I\cap F=\emptyset\).

    Therefore,
    \[
        \bigcap_{A\in F}d(A)\not\subseteq \bigcup_{B\in I}B.
    \]
    Hence there exists
    \[
        x\in
        \left(\bigcap_{A\in F}d(A)\right)
        \setminus
        \left(\bigcup_{B\in I}B\right).
    \]
    It follows that \(B\in I(x)\) for every \(B\in I\), and
    \(A\in F(x)\) for every \(A\in F\). Hence
    \[
        I\subseteq I(x)
        \qquad\text{and}\qquad
        F\subseteq F(x).
    \]
    Since \(I(x)\cap F(x)=\emptyset\), the maximality conditions in the
    definition of a comaximal pair imply that
    \[
        I=I(x)
        \qquad\text{and}\qquad
        F=F(x).
    \]

    \item[(v)] Finally, we prove that \((I(x),F(x))\) is a comaximal pair
    for every \(x\in X\), and that distinct points determine distinct
    pairs.

    By part (iii), \(I(x)\cap F(x)=\emptyset\). By the comaximal pair
    theorem, there exists a comaximal pair \((I,F)\) of
    \(\mathfrak{E}(X)\) such that
    \[
        I(x)\subseteq I
        \qquad\text{and}\qquad
        F(x)\subseteq F.
    \]
    By part (iv), there exists \(y\in X\) such that
    \[
        I=I(y)
        \qquad\text{and}\qquad
        F=F(y).
    \]
    Therefore
    \[
        I(x)\subseteq I(y)
        \qquad\text{and}\qquad
        F(x)\subseteq F(y).
    \]

    We claim that \(x=y\). The inclusion \(I(x)\subseteq I(y)\) implies
    that every essential subset containing \(y\) also contains \(x\).
    Since the essential subsets form a subbasis for \(\tau\), every
    \(\tau\)-neighborhood of \(y\) contains \(x\). Hence
    \[
        y\leq_{\tau}x.
    \]
    Similarly, the inclusion \(F(x)\subseteq F(y)\) implies that whenever
    \(x\in d(A)\), with \(A\) essential, then \(y\in d(A)\). Since the
    sets \(d(A)\), with \(A\) essential, form a basis for \(\sigma\),
    every \(\sigma\)-neighborhood of \(x\) contains \(y\). Hence
    \[
        x\leq_{\sigma}y.
    \]
    Since \(X\) is pairwise \(T_0\), the induced preorders are pairwise
    ordered. Therefore \(y\leq_{\tau}x\) and \(x\leq_{\sigma}y\) imply
    \(x=y\).

    Consequently,
    \[
        I=I(x)
        \qquad\text{and}\qquad
        F=F(x),
    \]
    so \((I(x),F(x))\) is a comaximal pair.

    Finally, if
    \[
        (I(x),F(x))=(I(y),F(y)),
    \]
    then the same argument applied to the inclusions
    \(I(x)\subseteq I(y)\) and \(F(x)\subseteq F(y)\) gives \(x=y\).
    Hence distinct points determine distinct comaximal pairs.
\end{enumerate}
\end{proof}
Define
\[
H_X:X\to spec_B(\mathfrak{E}(X)),\qquad x\mapsto (I(x),F(x)).
\]
The lemma shows $H_X$ is bijective.

\begin{Th}\label{HXHomeo}
For any pairwise Balbes-Dwinger space $(X,\tau,\sigma)$, the map $H_X$ is a bi-homeomorphism.
\end{Th}

\begin{proof}
By Lemma~\ref{CharMaxPair}, \(H_X\) is bijective. Let
\(A\in\mathfrak{E}(X)\). Then
\[
\begin{aligned}
x\in H_X^{-1}(\delta_{\mathfrak{E}(X)}(A))
&\Longleftrightarrow A\notin I(x)\\
&\Longleftrightarrow x\in A.
\end{aligned}
\]
Thus
\[
    H_X^{-1}(\delta_{\mathfrak{E}(X)}(A))=A.
\]
Since the essential subsets form a subbasis for \(\tau\), and the sets
\(\delta_{\mathfrak{E}(X)}(A)\) form a subbasis for the sup-topology of
\(spec_B(\mathfrak{E}(X))\), the map \(H_X\) is a homeomorphism with respect
to the first topologies.

Similarly,
\[
\begin{aligned}
x\in H_X^{-1}(\epsilon_{\mathfrak{E}(X)}(A))
&\Longleftrightarrow A\in F(x)\\
&\Longleftrightarrow x\in d_X(A).
\end{aligned}
\]
Hence
\[
    H_X^{-1}(\epsilon_{\mathfrak{E}(X)}(A))=d_X(A).
\]
Since the sets \(d_X(A)\), with \(A\in\mathfrak{E}(X)\), form a basis for
\(\sigma\), and the sets \(\epsilon_{\mathfrak{E}(X)}(A)\) form a basis for
the inf-topology of \(spec_B(\mathfrak{E}(X))\), the map \(H_X\) is also a
homeomorphism with respect to the second topologies. Therefore \(H_X\) is a
bi-homeomorphism.
\end{proof}

Hence any pairwise Balbes-Dwinger space is bihomeomorphic to the bitopological spectrum of a lattice.

\begin{Def}
A bicontinuous function
\[
    f:(X,\tau_1,\sigma_1)\to(Y,\tau_2,\sigma_2)
\]
between pairwise Balbes--Dwinger spaces is a \textbf{morphism of pairwise
Balbes--Dwinger spaces} if:
\begin{enumerate}
    \item[(i)] \(f\) is strongly bicontinuous;
    \item[(ii)] for every \(A\in\mathfrak{E}(Y)\),
    \[
        f^{-1}(d_Y(A))=d_X(f^{-1}(A));
    \]
    \item[(iii)] for every \(B\in d_Y(\mathfrak{E}(Y))\),
    \[
        f^{-1}(i_Y(B))=i_X(f^{-1}(B)).
    \]
\end{enumerate}
\end{Def}

\begin{Prop}\label{HXIso}
For any pairwise Balbes--Dwinger space \((X,\tau,\sigma)\), the map
\[
    H_X:X\rightarrow spec_B(\mathfrak{E}(X))
\]
is an isomorphism of pairwise Balbes--Dwinger spaces.
\end{Prop}

\begin{proof}
Let
\[
    S=spec_B(\mathfrak{E}(X)).
\]
For \(A\in\mathfrak{E}(X)\), write
\[
    \delta(A)=\delta_{\mathfrak{E}(X)}(A),
    \qquad
    \epsilon(A)=\epsilon_{\mathfrak{E}(X)}(A).
\]
By Theorem~\ref{HXHomeo}, the map \(H_X:X\to S\) is a bi-homeomorphism.
Hence \(H_X\) and \(H_X^{-1}\) are bicontinuous.

We first show that \(H_X\) is a morphism of pairwise Balbes--Dwinger spaces.
For every \(A\in\mathfrak{E}(X)\), we have
\[
\begin{aligned}
x\in H_X^{-1}(\delta(A))
&\Longleftrightarrow H_X(x)\in \delta(A)\\
&\Longleftrightarrow A\notin I(x)\\
&\Longleftrightarrow x\in A.
\end{aligned}
\]
Thus
\[
    H_X^{-1}(\delta(A))=A.
\]
Similarly,
\[
\begin{aligned}
x\in H_X^{-1}(\epsilon(A))
&\Longleftrightarrow H_X(x)\in \epsilon(A)\\
&\Longleftrightarrow A\in F(x)\\
&\Longleftrightarrow x\in d_X(A),
\end{aligned}
\]
and therefore
\[
    H_X^{-1}(\epsilon(A))=d_X(A).
\]

By Theorem~\ref{RepTh}, every essential subset of \(S\) is of the form
\(\delta(A)\), with \(A\in\mathfrak{E}(X)\). Since
\[
    H_X^{-1}(\delta(A))=A
\]
is essential in \(X\), the map \(H_X\) is strongly bicontinuous.

Now let \(A\in\mathfrak{E}(X)\). Since
\[
    d_S(\delta(A))=\epsilon(A),
\]
we have
\[
\begin{aligned}
H_X^{-1}(d_S(\delta(A)))
&=H_X^{-1}(\epsilon(A))\\
&=d_X(A)\\
&=d_X(H_X^{-1}(\delta(A))).
\end{aligned}
\]
Thus condition \((ii)\) in the definition of morphism holds for \(H_X\).

Next, every element of \(d_S(\mathfrak{E}(S))\) is of the form
\(\epsilon(A)\), with \(A\in\mathfrak{E}(X)\). For such an element we have
\[
    i_S(\epsilon(A))=\delta(A).
\]
Therefore
\[
\begin{aligned}
H_X^{-1}(i_S(\epsilon(A)))
&=H_X^{-1}(\delta(A))\\
&=A\\
&=i_X(d_X(A))\\
&=i_X(H_X^{-1}(\epsilon(A))).
\end{aligned}
\]
Thus condition \((iii)\) holds. Hence \(H_X\) is a morphism.

It remains to show that \(H_X^{-1}:S\to X\) is also a morphism. Since
\(H_X\) is bijective, for every \(A\in\mathfrak{E}(X)\) we have
\[
    (H_X^{-1})^{-1}(A)=H_X(A)=\delta(A).
\]
Indeed, \(H_X(x)\in\delta(A)\) if and only if \(x\in A\). Thus inverse
images of essential subsets of \(X\) under \(H_X^{-1}\) are essential subsets
of \(S\), and \(H_X^{-1}\) is strongly bicontinuous.

Let \(A\in\mathfrak{E}(X)\). Then
\[
\begin{aligned}
(H_X^{-1})^{-1}(d_X(A))
&=H_X(d_X(A))\\
&=\epsilon(A)\\
&=d_S(\delta(A))\\
&=d_S((H_X^{-1})^{-1}(A)).
\end{aligned}
\]
Thus condition \((ii)\) holds for \(H_X^{-1}\).

Finally, let \(B\in d_X(\mathfrak{E}(X))\). Then \(B=d_X(A)\) for some
\(A\in\mathfrak{E}(X)\). Since \(A\) is stable,
\[
    i_X(B)=i_X(d_X(A))=A.
\]
Hence
\[
\begin{aligned}
(H_X^{-1})^{-1}(i_X(B))
&=(H_X^{-1})^{-1}(A)\\
&=\delta(A).
\end{aligned}
\]
On the other hand,
\[
\begin{aligned}
i_S((H_X^{-1})^{-1}(B))
&=i_S((H_X^{-1})^{-1}(d_X(A)))\\
&=i_S(\epsilon(A))\\
&=\delta(A).
\end{aligned}
\]
Therefore
\[
    (H_X^{-1})^{-1}(i_X(B))
    =
    i_S((H_X^{-1})^{-1}(B)),
\]
so condition \((iii)\) also holds for \(H_X^{-1}\).

Consequently, both \(H_X\) and \(H_X^{-1}\) are morphisms of pairwise
Balbes--Dwinger spaces. Therefore \(H_X\) is an isomorphism in
\(\mathfrak{PBD}\).
\end{proof}
\subsection{Doubly Balbes-Dwinger spaces}
We now have a suitable category of spaces and morphisms. We first study what happens in the distributive case and characterize the pairwise Balbes-Dwinger spaces which come from a distributive lattice. Let us recall the function $b_L$, defined some sections above.
\[
b_L: spec(L)\rightarrow spec_B(L),\qquad P\mapsto (P,P^c).
\]

\begin{Prop}\label{topeq}
Let $L$ be a distributive lattice. Then:
\begin{enumerate}
    \item[(i)] $\tau_L=\sigma_L$, i.e., the sup-topology and the inf-topology of $spec_B(L)$ coincide.
    \item[(ii)] $b_L$ is a homeomorphism between $spec(L)$ and $(\mathfrak{M}(L),\tau_L)$.
\end{enumerate}
\end{Prop}

\begin{proof}
\begin{enumerate}
    \item[(i)] Since \(L\) is distributive, Proposition~\ref{TopDis} gives
    \[
        \delta_L=\epsilon_L.
    \]
    Hence
    \[
        \operatorname{Im}(\delta_L)=\operatorname{Im}(\epsilon_L).
    \]
    By definition, \(\tau_L\) is the topology generated by
    \(\operatorname{Im}(\delta_L)\) as a subbasis, while
    \(\operatorname{Im}(\epsilon_L)\) is a basis for \(\sigma_L\). Moreover,
    this common family is closed under finite intersections, since
    \(\operatorname{Im}(\epsilon_L)\) is a basis for \(\sigma_L\).  Therefore
    both topologies are generated by the same family of subsets of
    \(\mathfrak{M}(L)\), and hence
    \[
        \tau_L=\sigma_L.
    \]

    \item[(ii)] By Proposition~\ref{PCLat}, the map
    \[
        b_L:spec(L)\to\mathfrak{M}(L),\qquad P\mapsto (P,P^c),
    \]
    is bijective whenever \(L\) is distributive. Let \(x\in L\). Then
    \[
        P\in b_L^{-1}(\delta_L(x))
        \Longleftrightarrow
        (P,P^c)\in\delta_L(x)
        \Longleftrightarrow
        x\notin P
        \Longleftrightarrow
        P\in d_L(x).
    \]
    Thus
    \[
        b_L^{-1}(\delta_L(x))=d_L(x).
    \]
    Since the sets \(d_L(x)\) form a basis for the topology of \(spec(L)\),
    and the sets \(\delta_L(x)\) form a basis for \(\tau_L\) by part~\((i)\),
    it follows that \(b_L\) is continuous.

    Conversely, since \(b_L\) is bijective, the same computation gives
    \[
        b_L(d_L(x))=\delta_L(x)
    \]
    for every \(x\in L\). Hence \(b_L\) maps basic open subsets of
    \(spec(L)\) to basic open subsets of \((\mathfrak{M}(L),\tau_L)\).
    Therefore \(b_L\) is an open map. Consequently, \(b_L\) is a
    homeomorphism.
\end{enumerate}
\end{proof}

This proposition motivates the following definition.

\begin{Def}
A pairwise Balbes-Dwinger space $(X,\tau,\sigma)$ is a \textbf{doubly Balbes-Dwinger space} if $\tau=\sigma$.
\end{Def}

Then doubly Balbes-Dwinger spaces are the pairwise Balbes-Dwinger spaces which are essentially topological spaces, i.e., those for which one of the two topologies can be omitted.

\medskip

For a doubly Balbes--Dwinger space \((X,\tau,\sigma)\), the two topologies
coincide. In this case, the comparison with ordinary Balbes--Dwinger spaces
amounts to identifying the essential subsets with the corresponding
fundamental subsets of the common topology.

\begin{Prop}\label{BitoTop}
Let \((X,\tau,\sigma)\) be a pairwise Balbes--Dwinger space. If
\((X,\tau,\sigma)\) is doubly Balbes--Dwinger, then \((X,\tau)\) is a
Balbes--Dwinger space.
\end{Prop}

\begin{proof}
Since \(X\) is doubly Balbes--Dwinger, we have \(\tau=\sigma\). Hence the
two topological preorders coincide.

Let \(U\in\tau\). Since every open subset is increasing with respect to the
topological preorder, we have
\[
    \uparrow U=U.
\]
Moreover, \(U^c\) is decreasing. Hence
\[
    \downarrow(U^c)=U^c.
\]
Therefore
\[
    d(U)=\bigl(\downarrow(U^c)\bigr)^c=U
    \qquad\text{and}\qquad
    i(U)=\uparrow U=U.
\]

We now compare essential and fundamental subsets. Let \(A\) be a non-empty
essential subset. Then \(d(A)\) is open and \(A=i(d(A))\). Since \(i\) is the
identity on open subsets, we obtain
\[
    A=d(A).
\]
Thus \(A\) is open. Since \(A\) is also compact by definition, \(A\) is
fundamental.

Conversely, let \(A\) be a non-empty fundamental subset of \((X,\tau)\).
Then \(A\) is compact and open. Since \(d\) is the identity on open subsets,
we have \(d(A)=A\), which is open, and
\[
    i(d(A))=i(A)=A.
\]
Thus \(A\) is essential.

The definition of the empty essential subset also agrees with the definition
of the empty fundamental subset, because for open essential subsets we have
\(d(A)=A\). Hence the essential subsets of \(X\) coincide with the
fundamental subsets of \((X,\tau)\).

By the definition of pairwise Balbes--Dwinger space, essential subsets form a
subbasis for \(\tau\), and the family
\[
    \{d(A)\mid A\in\mathfrak{E}(X)\}
\]
is a basis for \(\sigma=\tau\), closed under finite intersections. Since
\(d(A)=A\) for fundamental subsets, the fundamental subsets form a basis for
\(\tau\), closed under finite intersections. Thus \((X,\tau)\) is coherent.

Moreover, \(d\)-birreducibility of essential subsets becomes ordinary
birreducibility of fundamental subsets, because \(d(A)=A\) for fundamental
\(A\). Finally, pairwise \(T_0\) reduces to the usual \(T_0\) condition when
\(\tau=\sigma\). Hence \((X,\tau)\) is a Balbes--Dwinger space.
\end{proof}

Thus doubly Balbes-Dwinger spaces correspond to Balbes-Dwinger spaces. It is natural, then, to study the behaviour of the lattice $\mathfrak{E}(X)$ for a doubly Balbes-Dwinger space $X$.

\begin{Prop}
If $(X,\tau,\sigma)$ is a doubly Balbes-Dwinger space, then $\mathfrak{E}(X)$ is a distributive lattice.
\end{Prop}

\begin{proof}
Let \(A,B\in\mathfrak{E}(X)\). Since \(X\) is doubly Balbes--Dwinger, the
two topologies coincide. By Proposition~\ref{BitoTop}, essential subsets
coincide with fundamental subsets of the common topology. Hence \(A\), \(B\),
and \(A\cap B\) are compact-open subsets of \((X,\tau)\).

Moreover, for open subsets we have \(d(U)=U\) and \(i(U)=U\). Therefore
\[
    A\wedge B=i(d(A\cap B))=A\cap B,
\]
while
\[
    A\vee B=A\cup B.
\]
Thus \(\mathfrak{E}(X)\) is a sublattice of \(\wp(X)\), and hence it is
distributive.
\end{proof}

The next theorem summarizes the behaviour of the distributive case in this theory.

\begin{Th}\label{DisChar}
Let $(X,\tau,\sigma)$ be a pairwise Balbes-Dwinger space. The following are equivalent:
\begin{enumerate}
    \item[(i)] $(X,\tau,\sigma)$ is a doubly Balbes-Dwinger space;
    \item[(ii)] $\mathfrak{E}(X)$ is a distributive lattice;
    \item[(iii)] There exists a distributive lattice $L$ such that $spec_B(L)$ is bi-homeomorphic to $X$;
    \item[(iv)] Every point of $(X,\tau,\sigma)$ is prime.
\end{enumerate}
\end{Th}
\begin{proof}
The implication \((i)\Rightarrow(ii)\) was proved in the preceding
proposition: if \(\tau=\sigma\), then essential subsets coincide with
fundamental subsets of the common topology, and the meet and join in
\(\mathfrak{E}(X)\) are intersection and union. Hence \(\mathfrak{E}(X)\) is
a sublattice of \(\wp(X)\).
If \(\mathfrak{E}(X)\) is distributive, then Theorem~\ref{HXHomeo} gives a
bi-homeomorphism
\[
    X \cong spec_B(\mathfrak{E}(X)),
\]
so \((ii)\Rightarrow(iii)\). Conversely, if \(X\) is bi-homeomorphic to
\(spec_B(L)\) for some distributive lattice \(L\), then
Proposition~\ref{topeq} gives \(\tau_L=\sigma_L\), and therefore
\(\tau=\sigma\). Hence \((iii)\Rightarrow(i)\).

It remains to relate these conditions with prime points. By
Theorem~\ref{HXHomeo}, \(X\) is bi-homeomorphic to
\(spec_B(\mathfrak{E}(X))\). Since being a prime point is invariant under
bi-homeomorphism, condition \((iv)\) holds for \(X\) if and only if every
point of \(spec_B(\mathfrak{E}(X))\) is prime. By
Theorem~\ref{PrimePointsDist}, this is equivalent to
\(\mathfrak{E}(X)\) being distributive. Thus \((iv)\Leftrightarrow(ii)\).
The result follows.
\end{proof}
\subsection{The main theorem}
In this section we will prove that $spec_B$ is a coequivalence of categories between the category of lattices with quasi-proper morphisms and the category of pairwise Balbes-Dwinger spaces. This is the main statement of this paper. We also review and unify the theory developed here with the work of Acosta, Balbes, Dwinger and Stone in the distributive case. As a consequence, we derive the duality theorem stated in Theorem~\ref{StoneDual}.

\begin{Prop}
If \(f:L\rightarrow N\) is a quasi-proper lattice homomorphism, then
\(spec_B(f):spec_B(N)\rightarrow spec_B(L)\) is a morphism of pairwise
Balbes--Dwinger spaces.
\end{Prop}

\begin{proof}
By Proposition~\ref{Func}, \(spec_B(f)\) is strongly bicontinuous.

Let \(A\in\mathfrak{E}(spec_B(L))\). By Theorem~\ref{RepTh}, there exists
\(x\in L\) such that
\[
    A=\delta_L(x).
\]
Then
\begin{align*}
spec_B(f)^{-1}(d_{spec_B(L)}(A))
&= spec_B(f)^{-1}(d_{spec_B(L)}(\delta_L(x)))\\
&= spec_B(f)^{-1}(\epsilon_L(x))\\
&= \epsilon_N(f(x))\\
&= d_{spec_B(N)}(\delta_N(f(x)))\\
&= d_{spec_B(N)}(spec_B(f)^{-1}(\delta_L(x)))\\
&= d_{spec_B(N)}(spec_B(f)^{-1}(A)).
\end{align*}
Thus condition \((ii)\) in the definition of morphism is satisfied.

Now let \(B\in d_{spec_B(L)}(\mathfrak{E}(spec_B(L)))\). Then
\(B=d_{spec_B(L)}(A)\) for some essential \(A\). Again, by
Theorem~\ref{RepTh}, there exists \(x\in L\) such that \(A=\delta_L(x)\).
Hence
\[
    B=d_{spec_B(L)}(\delta_L(x))=\epsilon_L(x).
\]
Therefore
\begin{align*}
spec_B(f)^{-1}(i_{spec_B(L)}(B))
&= spec_B(f)^{-1}(i_{spec_B(L)}(\epsilon_L(x)))\\
&= spec_B(f)^{-1}(\delta_L(x))\\
&= \delta_N(f(x))\\
&= i_{spec_B(N)}(\epsilon_N(f(x)))\\
&= i_{spec_B(N)}(spec_B(f)^{-1}(\epsilon_L(x)))\\
&= i_{spec_B(N)}(spec_B(f)^{-1}(B)).
\end{align*}
Thus condition \((iii)\) is satisfied.

Consequently, \(spec_B(f)\) is a morphism of pairwise Balbes--Dwinger spaces.
\end{proof}

Thus $spec_B$ can be regarded as a functor
\[
spec_B:\mathfrak{R}_{qp}\rightarrow \mathfrak{PBD}.
\]

On the other hand, for any morphism $f:X\rightarrow Y$ of pairwise Balbes-Dwinger spaces, define
\[
\mathfrak{E}(f)=f^{-1}:\mathfrak{E}(Y)\rightarrow\mathfrak{E}(X).
\]

\begin{Prop}
If \(f:X\rightarrow Y\) is a morphism of pairwise Balbes--Dwinger spaces, then
\[
    \mathfrak{E}(f)=f^{-1}:\mathfrak{E}(Y)\rightarrow\mathfrak{E}(X)
\]
is a quasi-proper lattice homomorphism.
\end{Prop}

\begin{proof}
First we show that \(\mathfrak{E}(f)\) is a lattice homomorphism. Let
\(A,B\in\mathfrak{E}(Y)\). Since \(f\) is strongly bicontinuous,
\(f^{-1}(A)\) and \(f^{-1}(B)\) are essential subsets of \(X\).

For joins, we have
\[
\begin{aligned}
\mathfrak{E}(f)(A\vee B)
&= f^{-1}(A\cup B)\\
&= f^{-1}(A)\cup f^{-1}(B)\\
&= \mathfrak{E}(f)(A)\vee \mathfrak{E}(f)(B).
\end{aligned}
\]

Now we prove preservation of meets. Since \(A\wedge B\) is essential and
stable, we have
\[
    A\wedge B=i_Y(d_Y(A\wedge B)).
\]
Therefore, using condition \((iii)\) in the definition of morphism with
\(d_Y(A\wedge B)\in d_Y(\mathfrak{E}(Y))\), we obtain
\[
\begin{aligned}
\mathfrak{E}(f)(A\wedge B)
&= f^{-1}(A\wedge B)\\
&= f^{-1}(i_Y(d_Y(A\wedge B)))\\
&= i_X(f^{-1}(d_Y(A\wedge B))).
\end{aligned}
\]
By Lemma~\ref{dMeet}, applied in \(Y\),
\[
    d_Y(A\wedge B)=d_Y(A)\cap d_Y(B).
\]
Hence
\[
\begin{aligned}
\mathfrak{E}(f)(A\wedge B)
&= i_X(f^{-1}(d_Y(A)\cap d_Y(B)))\\
&= i_X(f^{-1}(d_Y(A))\cap f^{-1}(d_Y(B)))\\
&= i_X(d_X(f^{-1}(A))\cap d_X(f^{-1}(B))),
\end{aligned}
\]
where the last equality follows from condition \((ii)\) in the definition of
morphism.

Applying Lemma~\ref{dMeet} in \(X\), we get
\[
    d_X(f^{-1}(A))\cap d_X(f^{-1}(B))
    =
    d_X(f^{-1}(A)\wedge f^{-1}(B)).
\]
Thus
\[
\begin{aligned}
\mathfrak{E}(f)(A\wedge B)
&= i_X(d_X(f^{-1}(A)\wedge f^{-1}(B)))\\
&= f^{-1}(A)\wedge f^{-1}(B),
\end{aligned}
\]
because \(f^{-1}(A)\wedge f^{-1}(B)\) is essential and hence stable. Therefore
\[
    \mathfrak{E}(f)(A\wedge B)
    =
    \mathfrak{E}(f)(A)\wedge \mathfrak{E}(f)(B).
\]
So \(\mathfrak{E}(f)\) is a lattice homomorphism.

It remains to show that \(\mathfrak{E}(f)\) is quasi-proper. Let \(x\in X\).
For every \(A\in\mathfrak{E}(Y)\), we have
\[
\begin{aligned}
A\in \mathfrak{E}(f)^{-1}(I(x))
&\Longleftrightarrow \mathfrak{E}(f)(A)\in I(x)\\
&\Longleftrightarrow x\notin f^{-1}(A)\\
&\Longleftrightarrow f(x)\notin A\\
&\Longleftrightarrow A\in I(f(x)).
\end{aligned}
\]
Hence
\[
    \mathfrak{E}(f)^{-1}(I(x))=I(f(x)).
\]

Similarly,
\[
\begin{aligned}
A\in \mathfrak{E}(f)^{-1}(F(x))
&\Longleftrightarrow \mathfrak{E}(f)(A)\in F(x)\\
&\Longleftrightarrow x\in d_X(f^{-1}(A))\\
&\Longleftrightarrow x\in f^{-1}(d_Y(A))\\
&\Longleftrightarrow f(x)\in d_Y(A)\\
&\Longleftrightarrow A\in F(f(x)),
\end{aligned}
\]
where we used condition \((ii)\) in the definition of morphism. Thus
\[
    \mathfrak{E}(f)^{-1}(F(x))=F(f(x)).
\]

By Lemma~\ref{CharMaxPair}, every comaximal pair of
\(\mathfrak{E}(X)\) is of the form \((I(x),F(x))\) for some \(x\in X\).
The equalities above show that its inverse image under \(\mathfrak{E}(f)\)
is
\[
    (I(f(x)),F(f(x))),
\]
which is a comaximal pair of \(\mathfrak{E}(Y)\), again by
Lemma~\ref{CharMaxPair}. Hence \(\mathfrak{E}(f)\) is quasi-proper.
\end{proof}

Thus $\mathfrak{E}$ can be regarded as a functor
\[
\mathfrak{E}:\mathfrak{PBD}\rightarrow\mathfrak{R}_{qp}.
\]

\begin{Prop}
The family $\{\delta_L:L\rightarrow \mathfrak{E}\circ spec_B(L)\}_{L\in Obj(\mathfrak{R}_{qp})}$ is a natural isomorphism between the identity functor of $\mathfrak{R}_{qp}$ and $\mathfrak{E}\circ spec_B$.
\end{Prop}

\begin{proof}
We already know that $\delta_L$ preserves joins. Let $x,y\in L$. Then
\begin{eqnarray*}
\delta_L(x)\wedge \delta_L(y)
&=& i(d(\delta_L(x)\cap \delta_L(y)))\\
&=& i(\epsilon_L(x)\cap\epsilon_L(y))\\
&=& i(\epsilon_L(x\wedge y))\\
&=& \delta_L(x\wedge y),
\end{eqnarray*}
so $\delta_L$ is a lattice homomorphism. By Theorem~\ref{RepTh}, $\delta_L$ is bijective. Hence $\delta_L$ is an isomorphism of lattices, and in particular a quasi-proper homomorphism.

It remains to prove naturality. Let \(f:L\rightarrow N\) be a quasi-proper
lattice homomorphism and let \(x\in L\). Then
\[
(\mathfrak{E}\circ spec_B)(f)(\delta_L(x))
=
\mathfrak{E}(spec_B(f))(\delta_L(x)).
\]
By definition, \(\mathfrak{E}(spec_B(f))\) is inverse image under
\(spec_B(f)\). Hence
\[
(\mathfrak{E}\circ spec_B)(f)(\delta_L(x))
=
spec_B(f)^{-1}(\delta_L(x)).
\]
By Proposition~\ref{Func}, we have
\[
spec_B(f)^{-1}(\delta_L(x))=\delta_N(f(x)).
\]
Therefore
\[
(\mathfrak{E}\circ spec_B)(f)(\delta_L(x))
=
\delta_N(f(x)).
\]
Thus
\[
(\mathfrak{E}\circ spec_B)(f)\circ\delta_L
=
\delta_N\circ f,
\]
i.e., the diagram
\[
\xymatrix{
L \ar[dd]_{f}\ar[rr]^{\delta_L} && \mathfrak{E}\circ spec_B(L)\ar[dd]^{\mathfrak{E}\circ spec_B(f)}\\
&&\\
N \ar[rr]_{\delta_N} && \mathfrak{E}\circ spec_B(N)
}
\]
commutes. This proves naturality.
\end{proof}

\begin{Prop}
The family
\[
    \{H_X:X\rightarrow spec_B(\mathfrak{E}(X))\}_{X\in Obj(\mathfrak{PBD})}
\]
is a natural isomorphism between the identity functor of \(\mathfrak{PBD}\)
and the functor \(spec_B\circ \mathfrak{E}\).
\end{Prop}

\begin{proof}
By Proposition~\ref{HXIso}, each \(H_X\) is an isomorphism of pairwise
Balbes--Dwinger spaces. It remains to prove naturality.

Let
\[
    f:X\rightarrow Y
\]
be a morphism of pairwise Balbes--Dwinger spaces. We must show that the
diagram
\[
\xymatrix{
X \ar[rr]^{H_X}\ar[dd]_{f} &&
spec_B(\mathfrak{E}(X)) \ar[dd]^{spec_B(\mathfrak{E}(f))}\\
&&\\
Y \ar[rr]_{H_Y} &&
spec_B(\mathfrak{E}(Y))
}
\]
commutes.

Let \(x\in X\). Write
\[
    H_X(x)=(I_X(x),F_X(x)),
    \qquad
    H_Y(f(x))=(I_Y(f(x)),F_Y(f(x))).
\]
Since
\[
    \mathfrak{E}(f)=f^{-1}:\mathfrak{E}(Y)\rightarrow\mathfrak{E}(X),
\]
we have
\[
spec_B(\mathfrak{E}(f))(H_X(x))
=
\left(
    \mathfrak{E}(f)^{-1}(I_X(x)),
    \mathfrak{E}(f)^{-1}(F_X(x))
\right).
\]

We compute the first component. For \(A\in\mathfrak{E}(Y)\),
\[
\begin{aligned}
A\in \mathfrak{E}(f)^{-1}(I_X(x))
&\Longleftrightarrow \mathfrak{E}(f)(A)\in I_X(x)\\
&\Longleftrightarrow f^{-1}(A)\in I_X(x)\\
&\Longleftrightarrow x\notin f^{-1}(A)\\
&\Longleftrightarrow f(x)\notin A\\
&\Longleftrightarrow A\in I_Y(f(x)).
\end{aligned}
\]
Hence
\[
    \mathfrak{E}(f)^{-1}(I_X(x))=I_Y(f(x)).
\]

Now we compute the second component. For \(A\in\mathfrak{E}(Y)\),
\[
\begin{aligned}
A\in \mathfrak{E}(f)^{-1}(F_X(x))
&\Longleftrightarrow \mathfrak{E}(f)(A)\in F_X(x)\\
&\Longleftrightarrow f^{-1}(A)\in F_X(x)\\
&\Longleftrightarrow x\in d_X(f^{-1}(A))\\
&\Longleftrightarrow x\in f^{-1}(d_Y(A))\\
&\Longleftrightarrow f(x)\in d_Y(A)\\
&\Longleftrightarrow A\in F_Y(f(x)),
\end{aligned}
\]
where we used condition \((ii)\) in the definition of morphism, namely
\[
    f^{-1}(d_Y(A))=d_X(f^{-1}(A)).
\]
Thus
\[
    \mathfrak{E}(f)^{-1}(F_X(x))=F_Y(f(x)).
\]

Therefore
\[
\begin{aligned}
spec_B(\mathfrak{E}(f))(H_X(x))
&=
\left(
    I_Y(f(x)),
    F_Y(f(x))
\right)\\
&=
H_Y(f(x)).
\end{aligned}
\]
Hence
\[
    spec_B(\mathfrak{E}(f))\circ H_X
    =
    H_Y\circ f.
\]
This proves naturality.
\end{proof}

We can now state the main theorem of this text, which summarizes a large part of the work developed here.

\begin{Th}[\textbf{Duality theorem for lattices}]\label{Dual}
The functor $spec_B:\mathfrak{R}_{qp}\rightarrow \mathfrak{PBD}$ is a coequivalence of categories, whose dual functor is $\mathfrak{E}:\mathfrak{PBD}\rightarrow\mathfrak{R}_{qp}$.
\end{Th}

We next examine how this theorem behaves in the distributive case, in order to verify that it is indeed a genuine generalization of Stone duality.

\medskip

Let $\mathfrak{DoBD}$ be the full subcategory of $\mathfrak{PBD}$ whose objects are doubly Balbes-Dwinger spaces, and let $\mathfrak{D}_p$ be the full subcategory of $\mathfrak{R}_{qp}$ whose objects are distributive lattices. Then $\mathfrak{D}_p$ is precisely the category of distributive lattices with proper homomorphisms. Theorem~\ref{Dual} and Theorem~\ref{DisChar} yield:

\begin{Th}[\textbf{Duality theorem for distributive lattices}]\label{DisDual}
The functor $spec_B$ restricts to a coequivalence
\[
spec_B:\mathfrak{D}_p\rightarrow \mathfrak{DoBD},
\]
whose dual functor is
\[
\mathfrak{E}:\mathfrak{DoBD}\rightarrow\mathfrak{D}_p.
\]
\end{Th}

Consider the functor
\[
\mathcal{O}:BiTop\rightarrow Top
\]
from the category of bitopological spaces with bicontinuous maps to the category of topological spaces, defined by $\mathcal{O}(X,\tau,\sigma)=(X,\tau)$ and $\mathcal{O}(f)=f$ on morphisms. Conversely, let
\[
\mathcal{D}:Top\rightarrow BiTop
\]
be the functor defined by $\mathcal{D}(X,\tau)=(X,\tau,\tau)$ and $\mathcal{D}(f)=f$ for continuous maps $f$.

By the definitions, if \((X,\tau)\) is a Balbes--Dwinger space, then
\((X,\tau,\tau)\) is a doubly Balbes--Dwinger space: in this case the
essential subsets coincide with the fundamental subsets of the common
topology. Conversely, Proposition~\ref{BitoTop} implies that if $(X,\tau,\sigma)$ is doubly Balbes-Dwinger, then $(X,\tau)$ is a Balbes-Dwinger space. Furthermore, if $f:X\rightarrow Y$ is a strongly continuous map between Balbes-Dwinger spaces, then $\mathcal{D}(f)$ is a morphism between doubly Balbes-Dwinger spaces; conversely, if $f:X\rightarrow Y$ is a morphism of doubly Balbes-Dwinger spaces, then $\mathcal{O}(f)$ is a strongly continuous map between Balbes-Dwinger spaces. Hence:

\begin{Prop}
$\mathcal{O}:\mathfrak{DoBD}\rightarrow\mathfrak{BD}$ is an isomorphism of categories, whose inverse is $\mathcal{D}:\mathfrak{BD}\rightarrow\mathfrak{DoBD}$.
\end{Prop}

Let $L$ be a distributive lattice. Recall
\[
b_L: spec(L)\rightarrow spec_B(L),\qquad P\mapsto (P,P^c),
\]
and note that we have already shown $b_L$ is a homeomorphism between $spec(L)$ and $\mathcal{O}\circ spec_B(L)$ (Proposition~\ref{topeq}). Moreover, it is straightforward to verify:

\begin{Prop}\label{NatIso}
The family $\{b_L:spec(L)\rightarrow \mathcal{O}\circ spec_B(L)\}_{L\in Obj(\mathfrak{D}_p)}$ is a natural isomorphism between the functor $spec$ and the functor $\mathcal{O}\circ spec_B$.
\end{Prop}

Therefore, Theorem~\ref{DisDual} and Proposition~\ref{NatIso} imply the following result.

\begin{Th}[\textbf{Stone duality theorem for distributive lattices}]
The functor $spec:\mathfrak{D}_p\rightarrow \mathfrak{BD}$ is a coequivalence of categories, whose dual functor is $\mathfrak{F}:\mathfrak{BD}\rightarrow\mathfrak{D}_p$.
\end{Th}

In other words, the Stone duality theorem (in the formulation of Acosta, Balbes and Dwinger) can be derived from the framework developed in this paper. Thus Theorem~\ref{Dual} is indeed a genuine generalization of Stone duality for distributive lattices to the setting of arbitrary lattices.

\section{Further directions}

The results of this paper extend the Stone--Balbes--Dwinger representation
and duality theorem from distributive lattices to arbitrary lattices. Several
natural directions remain to be explored.

\begin{Def}
A pairwise Balbes--Dwinger space \((X,\tau,\sigma)\) is \textbf{bounded} if
\((X,\tau)\) is compact and \(\emptyset\) is \(\sigma\)-fundamental.
\end{Def}

By combining Theorem~\ref{Dual} with Proposition~\ref{0IffFun} and
Proposition~\ref{1IffComp}, one obtains a coequivalence between bounded
pairwise Balbes--Dwinger spaces and bounded lattices with quasi-proper
homomorphisms. On the other hand, bounded lattices are represented by
Urquhart spaces \cite{Urquhart}, and the corresponding functorial version is
developed in \cite{Mancera}. Thus bounded pairwise Balbes--Dwinger spaces and
Urquhart spaces provide two different representations of the same algebraic
objects.

This suggests a more precise comparison between the two approaches. In the
distributive bounded case, Cornish \cite{Cornish} described the relationship
between Stone and Priestley representations. A natural continuation of the
present work is to establish an analogous comparison between bounded
pairwise Balbes--Dwinger spaces and Urquhart spaces. Such a comparison would
clarify how the bitopological Stone--Balbes--Dwinger viewpoint developed here
relates to the ordered Priestley--Urquhart viewpoint in the non-distributive
setting.

A second direction concerns functoriality. The duality theorem proved in this
paper represents all lattices as objects, but, as in the Balbes--Dwinger and
Urquhart settings, it requires a restriction on morphisms. The work of
Allwein and Hartonas \cite{AllHart} shows that, in the bounded case, full
functoriality can be recovered by considering all disjoint pairs of ideals
and filters rather than only comaximal pairs. It is therefore natural to ask
whether an analogous enlargement of the present bitopological construction
leads to a full duality for arbitrary lattices.

Another problem is to understand lattice-theoretic properties directly in
bitopological terms. In this paper, distributivity is characterized by the
coincidence of the two topologies, or equivalently by the collapse of the
bitopological spectrum to an ordinary Stone--Balbes--Dwinger space. It would
be interesting to find similar bitopological characterizations of other
important lattice-theoretic properties. For instance, one may ask for
necessary and sufficient conditions on a pairwise Balbes--Dwinger space
\(X\) ensuring that the lattice \(\mathfrak{E}(X)\) is modular.

Finally, the present construction shows that non-distributivity naturally
forces the Stone topology to split into two interacting topologies. This
raises a conceptual question: whether there exists a genuinely topological
representation of arbitrary lattices using only one topology, or whether the
bitopological structure is unavoidable for any Stone--Balbes--Dwinger type
extension beyond the distributive case.

\section*{Acknowledgements}
The author would like to express his gratitude to Professor Lorenzo Acosta G. for his valuable guidance, insightful discussions and encouragement during the development of this work.

This research received no external funding.

\end{document}